\documentclass[11pt]{elsarticle}

\usepackage[english]{babel}
\usepackage{graphicx,epstopdf,epsfig}
\usepackage{amsfonts,epsfig,fancyhdr,graphics, hyperref,amsmath,amssymb}
\usepackage{enumerate,mathtools,todonotes}
\hypersetup{colorlinks=true,
    linkcolor=blue,
    citecolor=blue,
    filecolor=cyan,      
    urlcolor=magenta
}
\usepackage[matrix,arrow,tips,curve]{xy}
\usepackage{latexsym,eucal,color}
\usepackage{pstricks}
\usepackage{dsfont}
\usepackage{setspace}
\usepackage{multicol}
\usepackage{tikz,pgfplots,graphicx}
\usetikzlibrary{patterns}
\usetikzlibrary{shapes,arrows,positioning,decorations.markings,math}
\usepackage{tkz-euclide}
\pgfplotsset{compat=1.11}
\usepackage[capitalise]{cleveref}

\newtheorem{theorem}{Theorem}[section]
\newtheorem{lemma}[theorem]{Lemma}
\newtheorem{prop}[theorem]{Proposition}
\newtheorem{corollary}[theorem]{Corollary}

\newtheorem{defn}[theorem]{Definition}

\newtheorem{remark}[theorem]{Remark}
\usepackage{subcaption}
\newcommand{\arcf}{\overrightarrow{\diamond}}
\newcommand{\arcb}{\overleftarrow{\diamond}}
\newcommand{\arcs}{\overleftrightarrow{\diamond}}
\newcommand{\corv}{\overleftrightarrow{\circ}}
\newcommand{\corvf}{\overrightarrow{\circ}}
\newcommand{\corvb}{\overleftarrow{\circ}}
\newcommand{\defin}[1]{\emph{#1}\index{#1}}

\usepackage{tikz}
\usetikzlibrary{arrows.meta,
	decorations.markings} 

\begin{document} 
\bibliographystyle{abbrv}
\setcounter{page}{1}

\thispagestyle{empty}

\title{Spectra of corona products of digraphs}

\date{~} 
  
\author[MC]{M. Cavers\texorpdfstring{\corref{cor1}}{}}
\ead{michael.cavers@utoronto.ca}
\author[FM]{F. Maghsoudi}
\ead{farzad.maghsoudi93@gmail.com}
\author[BM]{B. Miraftab}
\ead{bobby.miraftab@gmail.com}

\cortext[cor1]{Corresponding author}

\address[MC]{Department of Computer and Mathematical Sciences, University of Toronto Scarborough,\\ Toronto, Ontario, M1C 1A4, Canada}
\address[FM]{{Department of Mathematics and Computer Science, University of Lethbridge,\\ Lethbridge, Alberta, T1K 3M4, Canada}}
\address[BM]{School of Computer Science,
Carleton University,\\ Ottawa, Ontario, K1S 5B6, Canada}

\begin{abstract}
Two types of corona products for simple directed graphs are introduced, extending the classical notions from the undirected setting: the vertex-corona and the arc-corona. 
Their structural and spectral properties are investigated through the use of digraph coronals, with particular emphasis on the adjacency, Laplacian, and signless Laplacian spectra.
Finally, the coronals corresponding to these three matrices are computed for several families of digraphs.
\end{abstract}


\maketitle
\section{Introduction} \label{intro-sec}
In graph theory, various graph products, such as the disjoint union, Cartesian product, Kronecker (tensor) product, lexicographic product, direct product, and strong product, have been extensively studied. The spectral properties of these products, particularly concerning the adjacency, Laplacian, and signless Laplacian matrices, have been analyzed in numerous works \cite{MR3362696, cvetkovic1980spectra, godsil2001algebraic, imrich2000product}.
Among these, the corona product has gathered significant attention due to its unique structure and applications. Multiple variants of the corona product exist, and their adjacency, Laplacian, and signless Laplacian spectra have been explored in several studies 
\cite{BarikKalita, barik2007spectrum, barik2017laplacian, sharma2015structural, sonar2023duplication, sonar2023neighbourhood, sonar2024closed}.
However, the exploration of such graph products in the context of directed graphs (digraphs) is less known; existing results for digraph transformations appear in \cite{DengKelmans2013}, and for the Cartesian, lexicographic, direct and strong products in \cite{catral2020spectra}.
We refer the reader to \cite{Brualdi} for a survey on digraph spectra, and to \cite{BarikKalita} for one on adjacency, Laplacian, and signless Laplacian spectra of (undirected) graph coronas.
Notation and terminology not defined in this paper can be found in \cite{bang,Brualdi}.

The concept of the corona product of graphs was first introduced by Frucht and Harary \cite{frucht1970corona} in 1970, who defined the classical construction now known as the vertex corona. Since then, more than fifteen different variants of the corona product have been proposed in the literature, including edge coronas, neighbourhood coronas, subdivision-based coronas, and several double or composite coronas. This family of graph operations has attracted considerable attention from researchers, leading to a large body of work analyzing their structural and spectral properties, as well as numerous applications in graph theory and related fields.
In this paper, we generalize some of these definitions from graphs to digraphs, extending the theory to the directed setting, with the main focus on determining their spectra with respect to adjacency, Laplacian, and signless Laplacian matrices. 
Since a symmetric digraph corresponds to a graph, many of the results in this paper also hold for graphs.
\Cref{sec:prelim} covers terminology and basic results. 
\Cref{sec:coronals} examines the coronal of a digraph, its relation to the complement (\Cref{sec:comp}), and explicit formulas for families including joins of out-regular digraphs, semi-regular bipartite digraphs, directed paths 
and those whose adjacency matrices admit an equitable partition, i.e., a block matrix where each block has constant row sum (Section~\ref{sec:cor_comp}). 
Section~\ref{sec:vertexcor} defines the vertex corona and analyzes its characteristic polynomial. Section~\ref{sec:arccorona} introduces three arc-corona variants extending the edge-corona to digraphs. 

\subsection{Preliminaries}\label{sec:prelim}
Let $D$ be a digraph with vertex set $V(D) = \{v_1,v_2,\ldots,v_n\}$ and arc set $E(D)$ consisting of ordered pairs of vertices called \emph{arcs}. 
We use the notation $uv$ whenever $(u,v)\in E(D)$.
If both $uv,vu\in E(D)$, we call $uv$ a \emph{symmetric edge}. 
The \emph{underlying graph} $U(D)$ of $D$ satisfies $V(U(D))=V(D)$ and $E(U(D))=\{\{u,v\}:uv\in E(D)\}$. 
We say that $D$ is \emph{$r$-out-regular} if every vertex $v\in V(D)$ satisfies ${\rm deg}_{\rm out}(v)=r$, that is, every vertex has out-degree $r$.
The \emph{adjacency matrix} $A(D)=[a_{ij}]$ of $D$ is the $n \times n$ matrix with $a_{ij}=1$ if $v_iv_j\in E(D)$, and $a_{ij} = 0$ otherwise. 
The \emph{transpose} $\overleftarrow{D}$ of $D$ is obtained from $D$ by reversing the direction of each arc and has adjacency matrix $A(\overleftarrow{D})=A(D)^T$, where $M^T$ denotes the transpose of matrix $M$.
A \emph{symmetric digraph} is one whose adjacency matrix is symmetric (i.e., $A(D)=A(\overleftarrow{D})$) and can be identified with a simple graph.
Throughout this paper, we assume $D$ is a simple digraph (loops and multiple arcs are not permitted), that is, $A(D)$ is a $(0,1)$-matrix with zero diagonal.

We use $I_n$ to denote the $n\times n$ identity matrix, $O_{n_1,n_2}$ to denote the $n_1\times n_2$ zero matrix, 
$J_{n_1,n_2}$ to denote the $n_1\times n_2$ matrix with every element equal to one,
along with the shorthand $J_n$ to denote $J_{n,n}$, $O_n$ to denote $O_{n,n}$ and $\textbf{1}_n$ to denote $J_{n,1}$. 
The \emph{complement} $\overline{D}$ of $D$ is the simple digraph on the same vertex set $V(D)$ and $uv\in E(\overline{D})$ if and only if $uv\not\in E(D)$ for all $u,v\in V(D)$ with $v\neq u$. 
The adjacency matrix of the complement of $D$ is  $A(\overline{D})=J_n-A(D)-I_n$. 
The matrix $J_n-A(D)$ is the \emph{anti-adjacency matrix} of $D$.
The \emph{out-degree matrix} (resp. \emph{in-degree matrix}) $D_{\textnormal{out}}(D)$ (resp. $D_{\textnormal{in}}(D)$) is the diagonal matrix whose $i$-th diagonal entry equals the out-degree (resp. in-degree) of vertex $i$.
The \emph{Laplacian matrix} of $D$ is $L(D)=D_{\textnormal{out}}(D)-A(D)$
and the \emph{signless Laplacian matrix} of $D$ is $Q(D)=D_{\textnormal{out}}(D)+A(D)$. 
For an $n\times n$ matrix $M$, its \emph{characteristic polynomial} is defined as $f_M(\lambda) = \det(\lambda I_n - M)$. 
When $M$ is the adjacency (resp.\ Laplacian or signless Laplacian) matrix of a digraph $D$, the roots of $f_M(\lambda)$ are called the adjacency (resp.\ Laplacian or signless Laplacian) eigenvalues or spectrum of $D$. 

A digraph is \emph{bipartite} if its underlying graph is bipartite. 
A digraph is \emph{strongly connected} if it is possible to reach any vertex starting from any other vertex by traversing arcs in the direction(s) in which they point.
The line digraph of $D$, denoted by $\mathcal{L}(D)$, has one vertex for each arc of $D$, and two vertices representing arcs from $u$ to $v$ and from $w$ to $x$ in $D$ are connected by an arc from $uv$ to $wx$ in the line digraph if and only if $v = w$.
For a digraph $D$ with $V(D)=\{v_1,\ldots,v_n\}$ and arcs ordered as $e_1,\ldots,e_m$, its \emph{in-incidence matrix} $B_{\rm in}(D)=[b_{ij}]$ is an $n\times m$ matrix such that $b_{ij}=1$ if $e_j=v_kv_i$ for some vertex $v_k$, and $b_{ij}=0$ otherwise.
Similarly, the \emph{out-incidence matrix} of $D$ is the $n\times m$ matrix  $B_{\rm out}(D)=[b_{ij}']$ such that $b_{ij}'=1$ if $e_j=v_iv_k$ for some vertex $v_k$, and $b_{ij}'=0$ otherwise.
The product of the in- and out-incidence matrices of a digraph are related to the adjacency matrix of the digraph and its line digraph, as demonstrated in the next lemma.

 
\begin{lemma}{\rm\cite{Brualdi,MR2312328}}\label[lemma]{bb}
Let $D$ be a digraph and $\mathcal{L}(D)$ its line digraph. Then $A(D) = B_{\rm out}(D) B_{\rm in}(D)^T$ and $A(\mathcal{L}(D)) = B_{\rm in}(D)^T B_{\rm out}(D)$.
\end{lemma}

\cref{bb} mimics a result for graphs.
In particular, for a graph $G$ with $V(G)=\{v_1,\ldots,v_n\}$ and $E(G)=\{e_1,\ldots,e_m\}$, 
the \emph{incidence matrix} $B(G)=[b_{ij}]$ of $G$ is an $n\times m$ matrix such that $b_{ij}=1$ if $e_j=v_kv_i$ for some vertex $v_k$, and $b_{ij}=0$ otherwise. 
The \emph{oriented incidence matrix} of $G$, denoted by $N(G)$, is formed by first (arbitrarily) assigning an orientation to each edge of $G$ to produce a digraph $D$ (without $2$-cycles) and then computing $N(G)=B_{\rm out}(D)-B_{\rm in}(D)$.
It is known (e.g., see \cite{brouwer2011spectra}) that the Laplacian matrix of $G$ satisfies $L(G)=N(G)N(G)^T=D_{\rm deg}(G)-A(G)$ and the signless Laplacian matrix of $G$ satisfies $Q(G)=B(G)B(G)^T=D_{\rm deg}(G)+A(G)$, where $D_{\rm deg}(G)$ is the degree matrix of $G$ (i.e., the diagonal matrix whose $i$-th diagonal entry equals the degree of vertex $i$).
Furthermore, $A(\mathcal{L}(G))=B(G)^TB(G)-2I_m$, where $\mathcal{L}(G)$ is the line graph of $G$ (see \cite{cvetkovic1980spectra}).


In studying operations of graphs, techniques from linear algebra include properties of the Kronecker product, Schur complements (\cref{schur}), a consequence of Sylvester's determinant identity (\cref{sylv}(i)) and the Sherman-Morrison Formula (\cref{sylv}(ii)).
 
\begin{lemma}{\rm\cite[Eqn.~(6.2.1)]{meyer2023matrix}}
\label[lemma]{schur}
Let $ M_1, M_2, M_3 $, and $ M_4 $ be respectively $ p \times p $, $ p \times q $, $ q \times p $, and $ q \times q $ matrices.
If $M_4$ is invertible, then
\[\det \begin{bmatrix} 
    M_1 & M_2 \\
    M_3 & M_4 
\end{bmatrix} = \det(M_4) \cdot \det \left( M_1 - M_2 M_4^{-1} M_3 \right).\]
\end{lemma}
 
\begin{lemma}{\rm\cite[Eqn.~(6.2.3) and (3.8.2)]{meyer2023matrix}}\label[lemma]{sylv}
Let $C$ be an $n\times n$ invertible matrix and $\alpha\in\mathbb{R}$. Then 
\begin{enumerate}[(i)]
    \item 
$\displaystyle\det(C+\alpha J_n)=\det(C)\,(1+\alpha\,\textnormal{\textbf{1}}_n^TC^{-1}\textnormal{\textbf{1}}_n)$, 
\item\label{sylvii} 
$\displaystyle(C+\alpha J_n)^{-1}=C^{-1}-\frac{\alpha\,C^{-1}J_nC^{-1}}{1+\alpha \textnormal{\textbf{1}}_n^TC^{-1}\textnormal{\textbf{1}}_n}.$
\end{enumerate}
\end{lemma}

\section{Coronals of matrices and their computations}\label{sec:coronals}
Frucht and Harary \cite{frucht1970corona} first introduced the corona $G_1\circ G_2$ of two graphs $G_1$ and $G_2$.
We define and extend this concept to include digraphs in Section~\ref{sec:vertexcor}.
In order to describe the spectrum of $G_1 \circ G_2$,
McLeman and McNicholas \cite{Mcleman} introduced the coronal of a graph and demonstrated that the adjacency spectrum of $G_1 \circ G_2$ is completely determined by the spectra of $G_1$ and $G_2$ and the coronal of $G_2$.
They also computed the coronal of (undirected) paths, complete multipartite graphs, and joins of regular graphs. 
Since then, the notion of the coronal has been adapted to Laplacian matrices of graphs \cite{LIUu} and to real matrices \cite{cui2012spectrum}. 
We adopt the version for real matrices in this paper.
Specifically,
let $M$ be a real $n\times n$ matrix considered as a matrix over the field of rational functions $\mathbb{C}(\lambda)$ with $\det(\lambda I_n-M)$ nonzero.
The \emph{$M$-coronal} $\chi_M(\lambda) \in \mathbb{C}(\lambda)$ of $M$ is defined to be the sum of the entries of the matrix $(\lambda I_n-M)^{-1}$, i.e., $\chi_M(\lambda) = \textbf{1}_n^T\big(\lambda I_n-M\big)^{-1}\textbf{1}_n$. 

When $M$ is the adjacency (resp. Laplacian or signless Laplacian) matrix of a digraph $D$, we also refer to the $M$-coronal as the $A$-coronal (resp. $L$-coronal or $Q$-coronal) of $D$.
Note that a matrix and its transpose have the same coronal as well as any matrix permutation similar to it. 
The following lemma provides a simple closed-form expression for the coronal when the row sums of a matrix are constant but also applies if each column sum is constant.

\begin{lemma}{\rm\cite[Proposition 2]{cui2012spectrum}}\label[lemma]{conoal_cal}
Let $M$ be an $n\times n$ matrix such that each row sum of $M$ is equal to a constant $t$. Then $\chi_M (\lambda) = n/(\lambda-t)$.
\end{lemma}

Since the Laplacian matrix has constant row sum zero, \cref{conoal_cal} implies $\chi_{L(D)}(\lambda) = n/\lambda$ for any digraph $D$ with $n$ vertices.
This is identical to the graph case (see \cite{Liu}).  
In contrast, the $A$-coronal and $Q$-coronal are nontrivial. 



\subsection{Complements of digraphs via coronals}\label{sec:comp}
For a graph $G$ of order $n$, McLeman and McNicholas \cite[Theorem~12]{Mcleman} establish a relationship between the characteristic polynomials of $G$ and its complement $\overline{G}$ via the coronal of $G$, which can be rewritten as 
$\displaystyle f_{A(\overline{G})}(\lambda)=(-1)^nf_{A(G)}(-\lambda-1)\left(1+\chi_{A(G)}(-\lambda-1)\right)$.
To derive an analogous relationship for digraphs, we first prove an extension of this for real matrices. 
We also give a relationship between the coronal of an $n\times n$ matrix $M$ and of a linear combination of $M$, $I_n$ and $J_n$.

\begin{theorem}\label[theorem]{cor_addI}
Let $M$ be an $n\times n$ matrix and $M'=aM+bJ_n+cI_n$, where $a,b,c\in\mathbb{R}$ and $a\neq 0$.
Then the $M'$-coronal and the characteristic polynomial of $M'$ are, respectively, 
\begin{enumerate}[(i)]
    \item $\displaystyle\chi_{M'}(\lambda)= \frac{\chi_M\left(\tfrac{\lambda-c}{a}\right)}{a-b\, \chi_M\left(\tfrac{\lambda-c}{a}\right)}$,
    \item $\displaystyle f_{M'}(\lambda)= a^{n-1} f_M\left(\tfrac{\lambda-c}{a}\right)\left(a-b\, \chi_M\left(\tfrac{\lambda-c}{a}\right)\right).$
\end{enumerate}
\end{theorem}

\begin{proof}
For any $n\times n$ invertible matrix $C$ and constant $\alpha\in\mathbb{R}$, multiplying both sides of the equation in \cref{sylv}(ii) by $\textbf{1}_n^T$ (from the left) and $\textbf{1}_n$ (from the right) and noting $J_n=\textbf{1}_n\textbf{1}_n^T$ gives
\begin{align}
\textbf{1}_n^T(C+\alpha J_n)^{-1}\textbf{1}_n
&=\textbf{1}_n^TC^{-1}\textbf{1}_n-\frac{\alpha\,\textbf{1}_n^TC^{-1}\textbf{1}_n\textbf{1}_n^TC^{-1}\textbf{1}_n^T}{1+\alpha \textbf{1}_n^TC^{-1}\textbf{1}_n}
=\frac{\textbf{1}_n^TC^{-1}\textbf{1}_n}{1+\alpha \textbf{1}_n^TC^{-1}\textbf{1}_n}.\label{eq1}
\end{align}
Let $\alpha=-b$ and $C=(\lambda-c)I_n-aM$ (and consider $C$ as an invertible matrix over the field of rational functions $\mathbb{C}(\lambda)$).
Then the left-side of \eqref{eq1} is
\begin{align*}
\textbf{1}_n^T(C+\alpha J_n)^{-1}\textbf{1}_n
&=\textbf{1}_n^T((\lambda-c)I_n-aM-bJ_n)^{-1}\textbf{1}_n\\
&=\textbf{1}_n^T(\lambda I_n-(aM+bJ_n+cI_n))^{-1}\textbf{1}_n\\
&=\chi_{M'}(\lambda).
\end{align*}
Next observe that
\begin{align}
\textbf{1}_n^TC^{-1}\textbf{1}_n
&=\textbf{1}_n^T\left((\lambda-c)I_n-aM\right)^{-1}\textbf{1}_n\label{eq2}\\ 
&=\frac{1}{a}\textbf{1}_n^T\left(\tfrac{\lambda-c}{a}\right)I_n-M)^{-1}\textbf{1}_n\nonumber\\ 
&=\frac{1}{a}\chi_M\left(\tfrac{\lambda-c}{a}\right).\nonumber
\end{align}
The result now follows by substituting \eqref{eq2} into the right-side of \eqref{eq1} and multiplying both the numerator and denominator by $a$.
For the second equation, by \cref{sylv}(i), we obtain
\begin{align*}
f_{M'}(a\lambda+c) &= \det\left((a\lambda+c)I_n-(aM+bJ_n+cI_n)\right)\\
&=a^n\det\left((\lambda I_n-M)-(b/a)J_n\right)\\
&=a^n\det(\lambda I_n-M)\,\left(1-(b/a)\textbf{1}_n^T(\lambda I_n-M)^{-1}\textbf{1}_n\right)\\
&=a^n\,f_M(\lambda)\,\left(1-(b/a)\chi_M(\lambda)\right)\\
&=a^{n-1}\,f_M(\lambda)\,\left(a-b\,\chi_M(\lambda)\right).
\end{align*} 
The result follows by replacing $\lambda$ by $\tfrac{\lambda-c}{a}$.
\end{proof}



A special case of \cref{cor_addI}(ii) is given by Liu and Zhang \cite[Corollary~2.3]{Liu} (with $a=1$, $b=\alpha$ and $c=0$).
We next obtain a relationship between the $A$-coronal (resp. $Q$-coronal) of a digraph and that of its complement. 

\begin{corollary}\label[corollary]{cor-comp}
Let $D$ be a digraph with $n$ vertices and $\overline{D}$ its complement. Then
\begin{enumerate}[(i)]
    \item $\displaystyle\chi_{A(\overline{D})}(\lambda)=\frac{1}{1+\chi_{A(D)}(-\lambda-1)}-1$,
    \item $\displaystyle\chi_{Q(\overline{D})}(\lambda)=\frac{1}{1+\chi_{Q(D)}(n-\lambda-2)}-1$.
\end{enumerate}
\end{corollary}
\begin{proof}
For (i), apply \cref{cor_addI}(i) where $M=A(D)$ and $M'=-M+J_n-I_n$, and note that $M'=A(\overline{D})$, and for (ii), take $M=Q(D)$ and $M'=-M+J_n+(n-2)I_n$, and note that $M'=Q(\overline{D})$. 
\end{proof}

Formulas for the adjacency (resp. Laplacian and signless Laplacian) characteristic polynomial of the complement of a digraph in terms of the respective coronals can also be obtained using \cref{cor_addI}.

\begin{corollary}\label[corollary]{prop:cor_comp}
Let $D$ be a digraph with $n$ vertices and $\overline{D}$ its complement. Then
\begin{enumerate}[(i)]
    \item $\displaystyle f_{A(\overline{D})}(\lambda) = (-1)^n\big(1 + \chi_{A(D)}(-\lambda - 1)\big)\, f_{A(D)}(-\lambda - 1)$,
    \item $\displaystyle f_{L(\overline{D})}(\lambda)=(-1)^n\frac{\lambda}{\lambda-n}f_{L(D)}(n-\lambda)$,
    \item $\displaystyle f_{Q(\overline{D})}(\lambda) = (-1)^n\big(1 + \chi_{Q(D)}(n-\lambda-2)\big)\, f_{Q(D)}(n-\lambda-2)$.
\end{enumerate}
\end{corollary}

\begin{proof}
In each case, apply \cref{cor_addI}(ii).
For the first formula, take $M=A(D)$ and $M'=-M+J_n-I_n$, and note that $M'=A(\overline{D})$. 
For the second formula, take $M=L(D)$ and $M'=-M-J_n+nI_n$, and note that $M'=L(\overline{D})$ and $\chi_{L(D)}(\lambda) = n/\lambda$ (by \cref{conoal_cal}).
For the third formula, take $M=Q(D)$ and $M'=-M+J_n+(n-2)I_n$, and note that $M'=Q(\overline{D})$. 
\end{proof}

The characteristic polynomial of the Laplacian of the complement of a digraph $D$ (\cref{prop:cor_comp}(ii)) is also stated in  
\cite[Theorem~3.6]{deng2017laplacian} and \cite{Kelmans1965}, and also proven in \cite[Eq.~(6)]{su2012laplacian} using an alternate argument.
When $D$ has a simple $A$-coronal or $Q$-coronal, \cref{prop:cor_comp} yields a simpler formula, for example, in the case that $D$ is out-regular.
Note that \cref{cor:charcomp}(i) also appears in \cite[Corollary~3.12]{DengKelmans2013}.

\begin{corollary}\label[corollary]{cor:charcomp}
Let $D$ be an $r$-out-regular digraph with $n$ vertices and let $\overline{D}$ be its complement. Then
\begin{enumerate}[(i)]
    \item $\displaystyle f_{A(\overline{D})}(\lambda)= (-1)^n \frac{\lambda - n + r + 1}{\lambda + r + 1} f_{A(D)}(-\lambda - 1)$,
    \item $\displaystyle f_{Q(\overline{D})}(\lambda)= (-1)^n \frac{\lambda -2n + 2r + 2}{\lambda -n + 2r + 2} f_{Q(D)}(n-\lambda-2)$.
\end{enumerate}
\end{corollary}

\begin{proof}
Apply \cref{prop:cor_comp} and note that (by \cref{conoal_cal}) $\chi_{A(D)}(\lambda)=n/(\lambda-r)$ and $\chi_{Q(D)}(\lambda)=n/(\lambda-2r)$ since $D$ is $r$-out-regular.
\end{proof}

\subsection{Computations of coronals}\label[corollary]{sec:cor_comp}
Since coronals play a key role in the study of digraph coronas and digraph complements, we compute the $A$-coronal and $Q$-coronal for several digraph families in this section.  
We omit the $L$-coronal as \cref{conoal_cal} implies $\chi_{L(D)}(\lambda) = n/\lambda$ for any digraph $D$ with $n$ vertices.
As mentioned previously, when $D$ is $r$-out-regular, \cref{conoal_cal} gives $\chi_{A(D)}(\lambda)=n/(\lambda-r)$ and $\chi_{Q(D)}(\lambda)=n/(\lambda-2r)$.
For example, the directed cycle $C_n$ with $n$ vertices satisfies
$\chi_{A(C_n)}(\lambda)=n/(\lambda-1)$
and $\chi_{Q(C_n)}(\lambda)=n/(\lambda-2)$.

\cref{conoal_cal} admits a partial extension to matrices with equitable partitions.
Specifically, a (block) partition of a square matrix $M$ is called \emph{equitable} if each block has constant row sums and all diagonal blocks are square. The corresponding \emph{quotient matrix} $R$ is defined by taking, as its $(i,j)$-entry, the common row sum of the block in position $(i,j)$ of $M$.

\begin{lemma}\label[lemma]{thm:equitable}
Let $k\geq 1$ and 
\[
M=\begin{bmatrix} 
    M_{1,1} & \cdots & M_{1,k} \\
    \vdots & & \vdots\\
    M_{k,1} & \cdots & M_{k,k}
\end{bmatrix},
\]
where each block $M_{i,j}$ is an $n_i\times n_j$ matrix with constant row sum $r_{ij}$ for $1\leq i,j\leq k$. 
Suppose $R=[r_{ij}]$ is the $k\times k$ quotient matrix of $M$.
Then
\[
\chi_M(\lambda)=\begin{bmatrix}
    n_1&n_2&\cdots&n_k
\end{bmatrix}
\left(\lambda I_k-R\right)^{-1}\textnormal{\textbf{1}}_k.
\]
In particular, if $k=2$ then 
\[
\chi_M(\lambda)=\frac{(n_1+n_2)\lambda+n_1(r_{12}-r_{22})+n_2(r_{21}-r_{11})}{\lambda^2-(r_{11}+r_{22})\lambda+r_{11}r_{22}-r_{12}r_{21}}.
\]
\end{lemma}

\begin{proof}
Let $n=n_1+\cdots+n_k$ and $X={\rm diag}(a_1I_{n_1},\ldots,a_kI_{n_k})$, where $\begin{bmatrix}
    a_1&
    \cdots&
    a_k
\end{bmatrix}^T=\left(\lambda I_k-R\right)^{-1}\textbf{1}_{k}$.
This means 
$\left(\lambda I_k-R\right)\begin{bmatrix}
    a_1&
    \cdots&
    a_k
\end{bmatrix}^T=\textbf{1}_{k}$, and so, $\lambda a_i-\sum_{j=1}^ka_jr_{ij}=1$ for all $1\leq i\leq k$.
Combining this with $M_{i,j}\textbf{1}_{n_j}=r_{ij}\textbf{1}_{n_j}$,
we obtain $(\lambda I_n - M) X \textbf{1}_{n}=\textbf{1}_{n}$.
Thus,
\[
\chi_M(\lambda)=\textbf{1}_{k}^TX\textbf{1}_{k}=n_1a_1+n_2a_2+\cdots+n_ka_k
=\begin{bmatrix}
    n_1&n_2&\cdots&n_k
\end{bmatrix}
\left(\lambda I_k-R\right)^{-1}\textbf{1}_{k}.
\]
The formula when $k=2$ follows directly from the inverse formula for a $2\times 2$ matrix.
\end{proof}

The case $k=2$ in \cref{thm:equitable} also appears in \cite[Theorem~2.7]{RajkumarGayathri2020} with an alternate proof.
\cref{cor_addI} and \cref{thm:equitable} can be used to 
derive a formula for the join of out-regular digraphs, giving a generalization to the formula computed in \cite[Proposition~17]{Mcleman} for the $A$-coronal of the join of two regular graphs.
Following \cite{cvetkovic1980spectra}, we define the join (also called the complete product) of two digraphs $D_1$ and $D_2$, denoted $D_1\vee D_2$, as the digraph obtained from their disjoint union by adding all possible arcs between every vertex of $D_1$ and every vertex of $D_2$.

\begin{prop}\label[prop]{prop:join}
Let $k\geq 1$. For each $i=1,\ldots, k$, let $D_i$ be an $r_i$-out-regular digraph with $n_i$ vertices.
Then
\begin{enumerate}[(i)]
    \item $\displaystyle \chi_{A(D_1\vee\cdots\vee D_k)}(\lambda)=\left[1-\sum_{i=1}^k\frac{n_i}{\lambda+n_i-r_i}\right]^{-1}-1$,
    \item $\displaystyle \chi_{Q(D_1\vee\cdots\vee D_k)}(\lambda)=\left[1-\sum_{i=1}^k\frac{n_i}{\lambda-n+2n_i-2r_i}\right]^{-1}-1$, where $n=n_1+\cdots+n_k$.
\end{enumerate}

\end{prop}

\begin{proof}
Let $D=D_1\vee\cdots\vee D_k$.
We first compute $\chi_{A(\overline{D})}(\lambda)$.
Since $\overline{D}$ is the disjoint union of $\overline{D_1},\ldots,\overline{D_k}$, its adjacency matrix is permutation similar to 
$A(\overline{D_1})\oplus\cdots\oplus A(\overline{D_k})$, where $M_1\oplus M_2$ denotes the direct sum of square matrices $M_1$ and $M_2$.
This matrix admits an equitable partition with quotient matrix
$R_A={\rm diag}(n_1-r_1-1,\ldots,n_k-r_k-1)$.
By \cref{thm:equitable}, 
\begin{align*}
\chi_{A(\overline{D})}(\lambda)
&=\begin{bmatrix}
    n_1&n_2&\cdots&n_k
\end{bmatrix}
\left(\lambda I_k-R_A\right)^{-1}\textbf{1}_k\\
&=\begin{bmatrix}
    n_1&n_2&\cdots&n_k
\end{bmatrix}
{\rm diag}(1/(\lambda-n_1+r_1+1),\ldots,1/(\lambda-n_k+r_k+1))\textbf{1}_k\\
&=\sum_{i=1}^k\frac{n_i}{\lambda-n_i+r_i+1}.
\end{align*}
By \cref{cor-comp}(i),
\[
\chi_{A(D)}(\lambda)
=\left[1+\sum_{i=1}^k\frac{n_i}{(-\lambda-1)-n_i+r_i+1}\right]^{-1}-1
=\left[1-\sum_{i=1}^k\frac{n_i}{\lambda+n_i-r_i}\right]^{-1}-1.
\]
The argument for $\chi_{Q(\overline{D})}(\lambda)$ is identical:
$Q(\overline{D})$ is permutation similar to 
$Q(\overline{D_1})\oplus\cdots\oplus Q(\overline{D_k})$, which admits an equitable partition with quotient matrix
$R_Q=2\cdot{\rm diag}(n_1-r_1-1,\ldots,n_k-r_k-1)$.
By \cref{thm:equitable}, 
\begin{align*}
\chi_{Q(\overline{D})}(\lambda)
&=\begin{bmatrix}
    n_1&n_2&\cdots&n_k
\end{bmatrix}
\left(\lambda I_k-R_Q\right)^{-1}\textbf{1}_k\\
&=\begin{bmatrix}
    n_1&n_2&\cdots&n_k
\end{bmatrix}
{\rm diag}(1/(\lambda-2n_1+2r_1+2),\ldots,1/(\lambda-2n_k+2r_k+2))\textbf{1}_k\\
&=\sum_{i=1}^k\frac{n_i}{\lambda-2n_i+2r_i+2}.
\end{align*}
By \cref{cor-comp}(ii),
\[
\chi_{Q(D)}(\lambda)
=\left[1+\sum_{i=1}^k\frac{n_i}{(n-\lambda-2)-2n_i+2r_i+2}\right]^{-1}-1
=\left[1-\sum_{i=1}^k\frac{n_i}{\lambda-n+2n_i-2r_i}\right]^{-1}-1.
\]
\end{proof}

\cref{thm:equitable} can also be applied to other classes of digraphs such as \emph{semi-regular} bipartite digraphs, that is, 
bipartite digraphs whose bipartition $(V_1,V_2)$ has the property that every vertex in $V_1$ has constant out-degree and every vertex in $V_2$ has constant out-degree (these two constants can be different).

\begin{prop} 
Let $D$ be a semi-regular bipartite digraph with bipartition $(V_1,V_2)$ where every vertex in $V_1$ has out-degree $r_1$ and every vertex in $V_2$ has out-degree $r_2$. Suppose $n_1=|V_1|$ and $n_2=|V_2|$. Then
\begin{enumerate}[(i)]
    \item $\displaystyle \chi_{A(D)}(\lambda)=\frac{(n_1+n_2)\lambda+n_1r_{1}+n_2r_{2}}{\lambda^2-r_{1}r_{2}}$,
    \item $\displaystyle \chi_{Q(D)}(\lambda)=\frac{(n_1+n_2)\lambda+(n_1-n_2)(r_{1}-r_{2})}{\lambda(\lambda-(r_{1}+r_{2}))}$.
\end{enumerate}
\end{prop}
 
\begin{proof}
The adjacency and signless Laplacian matrices for $D$ each admit an equitable partition with (up to permutation) quotient matrices 
$
R_A=
\left[\begin{array}{cc}
0 & r_1 \\
r_2 & 0\\
\end{array}\right]$ and 
$R_Q=
\left[\begin{array}{cc}
r_1 & r_1 \\
r_2 & r_2\\
\end{array}\right]$,
respectively.
The result now follows from \cref{thm:equitable}.
\end{proof}

In cases where formulas are known for both the characteristic polynomial of $D$ and either its complement or anti-adjacency matrix, \cref{cor_addI} can be applied to provide a simple formula for its coronal.
We demonstrate with another class of bipartite digraphs that admit a simple $A$-coronal. 

\begin{prop} 
Let $D$ be a bipartite digraph with bipartition $(V_1, V_2)$ and suppose $n_1=|V_1|$, $n_2=|V_2|$.
If 
${\rm deg}_{\rm out}(v)=n_2$ for every $v\in V_1$ and 
$k=\sum_{v\in V_2}{\rm deg}_{\rm out}(v)$,
then
\[
\chi_{A(D)}(\lambda) = \frac{(n_1+n_2)\lambda + k + n_1n_2}{\lambda^2 - k}.
\]
\end{prop}

\begin{proof}
Up to permutation similarity, the adjacency matrix of $D$ has the form
$
A(D) = \begin{bmatrix} 
    0 & J_{n_1,n_2} \\
    B & 0 
\end{bmatrix},
$
where 
$B$ is a $(0,1)$-matrix of size $n_2 \times n_1$ with exactly $k$ ones. 
By \cite[Lemma~4.1]{cavers2025digraphs}, the characteristic polynomial of $A(D)$ is
$f_{A(D)}(\lambda) = \lambda^{n_1+n_2-2}(\lambda^2 - k)$.
Since the anti-adjacency matrix $C=J_{n_1+n_2}-A(D)$ of $D$ is lower block triangular with diagonal blocks $J_{n_1}$ and $J_{n_2}$, its characteristic polynomial is
$f_{C}(\lambda) = \lambda^{n_1+n_2-2}(\lambda - n_1)(\lambda - n_2)$.
Applying \cref{cor_addI} with $a = -1$, $b = 1$, and $c = 0$, we obtain
\[
\chi_{A(D)}(\lambda)
= -1 + \frac{(-\lambda)^{n_1+n_2-2}(-\lambda - n_1)(-\lambda - n_2)}{(-1)^{n_1+n_2} \lambda^{n_1+n_2-2} (\lambda^2 - k)}
= \frac{(n_1+n_2)\lambda + k + n_1n_2}{\lambda^2 - k}.
\]
\end{proof}

For the directed path, closed-form formulas for the $A$-coronal and $Q$-coronal can be derived directly from the definition.  

\begin{prop}\label[prop]{prop:coronoal_path}
Let $P_n$ denote the directed path on $n$ vertices.
Then 
\begin{enumerate}[(i)]
    \item $\displaystyle \chi_{A(P_n)}(\lambda) = \frac{n\lambda^{n+1}-(n+1)\lambda^n+1}{\lambda^n(\lambda-1)^2}$,
    \item $\displaystyle \chi_{Q(P_n)}(\lambda)=\frac{(\lambda-1)^n\big(\lambda n(\lambda-2)-2(\lambda-1)\big)+2(\lambda-1)}{\lambda(\lambda-1)^n(\lambda-2)^2}$.
\end{enumerate}
\end{prop}
 
\begin{proof}
In \cite{usmani1994inversion}, Usmani provided an explicit formula for the inverse of a tridiagonal matrix from which the inverse of a bidiagonal matrix can be computed.
It follows that if $B=\lambda I_n-A(P_n)$, then 
$(B^{-1})_{ij}=\lambda^{i-j-1}$ if $1\leq i\leq j\leq n$, and $(B^{-1})_{ij}=0$ otherwise.
Thus, 
\[
\chi_{A(P_n)}(\lambda) = \sum_{i=1}^n\sum_{j=i}^n
\lambda^{i-j-1}
=\sum_{i=1}^n
\lambda^{i-1}\frac{\lambda^{-i}(1-\lambda^{-(n-i+1)})}{1-\lambda^{-1}} 
=\frac{n\lambda^{-1}}{1-\lambda^{-1}}-
\frac{\lambda^{-1}\lambda^{-n-1}}{1-\lambda^{-1}}\sum_{i=1}^n\lambda^i,
\]
from which the result for the $A$-coronal of $P_n$ follows (by applying the formula for the sum of a geometric sequence).
For the signless Laplacian, if $B=\lambda I_n-Q(P_n)$, then 
$(B^{-1})_{ij}=(\lambda-1)^{i-j-1}$ if $1\leq i\leq j\leq n-1$, 
$(B^{-1})_{i,n}=\lambda^{-1}(\lambda-1)^{i-n}$ if $1\leq i\leq n$, 
and $(B^{-1})_{ij}=0$ otherwise.
Thus, 
\[
\chi_{Q(P_n)}(\lambda) 
= \sum_{i=1}^{n-1}\sum_{j=i}^{n-1}
(\lambda-1)^{i-j-1}
+\sum_{i=1}^n\lambda^{-1}(\lambda-1)^{i-n},
\]
from which the result for the $Q$-coronal of $P_n$ follows (again by applying the formula for the sum of a geometric sequence multiple times).
\end{proof}

The first equation in \cref{prop:coronoal_path} may be rewritten as $\chi_{A(P_n)}(\lambda)=\left(\sum_{k=1}^{n}k\lambda^{k}\right)/\lambda^{n+1}$
and can also be proved using \cref{cor_addI} and 
\cite[Lemma~4.2.1]{prayitno2025eigenvalues} (the latter result gives a formula for the characteristic polynomial of the anti-adjacency matrix of the directed path $P_n$).

\section{Vertex Corona}\label{sec:vertexcor}
Frucht and Harary \cite{frucht1970corona} first introduced the corona of two graphs:
Given two graphs $G_1$ and $G_2$ where $G_1$ has $n_1$ vertices, the corona $G_1 \circ G_2$ is the graph obtained by taking one copy of $G_1$ and $n_1$ copies of $G_2$, and then joining the $i$-th vertex of $G_1$ to every vertex in the $i$-th copy of $G_2$.
We propose the following extensions to digraphs. 

\begin{defn}\rm \label{2-cycle-vertex}
Let $D_1$ and $D_2$ be digraphs with $D_1$ having $n_1$ vertices, say $v_1,\ldots,v_{n_1}$. 
The \defin{forward-vertex-corona} $D_1 \corvf D_2$ is formed by taking one copy of $D_1$ and $n_1$ copies of $D_2$, and for each vertex $v_k$ in $V(D_1)$ and each vertex $w$ in the $k$-th copy of $D_2$, adding the arc $v_kw$.
The \defin{backward-vertex-corona} $D_1 \corvb D_2$ instead adds the arc $wv_k$ while the \defin{symmetric-vertex-corona} $D_1 \corv D_2$ adds both $v_kw$ and $wv_k$. 
\end{defn}

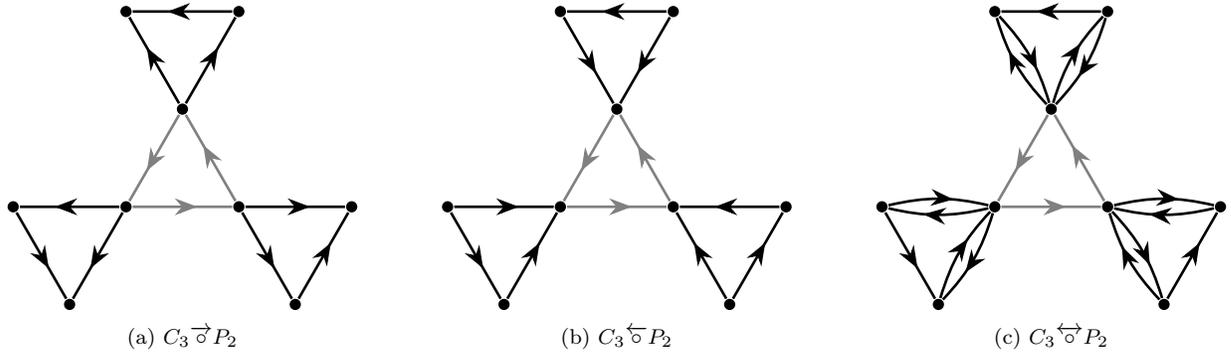
\begin{figure}[!ht]
    \centering  
    \begin{subfigure}[b]{0.3\textwidth}
        \centering
                \begin{tikzpicture}[decoration = {markings,mark=at position .65 with {\arrow{Stealth[length=3mm]}}},
            scale=0.75,
            v/.style={circle,fill,inner sep=1.5pt},
            thick,every edge/.style = {draw, line width=1.pt, postaction=decorate},
        ]
          \pgfmathsetmacro{\hs}{sqrt(3)}  
          
          \node[v] (T1) at (-1,  3) {};            
          \node[v] (T2) at ( 1,  3) {};            
          \node[v] (C)  at ( 0, {3 - \hs}) {};
          
          \node[v] (M1) at (-3, {3 - 2*\hs}) {};   
          \node[v] (M2) at (-1, {3 - 2*\hs}) {};   
          \node[v] (M3) at ( 1, {3 - 2*\hs}) {};   
          \node[v] (M4) at ( 3, {3 - 2*\hs}) {};   
          
          \node[v] (BL) at (-2, {3 - 3*\hs}) {};  
          \node[v] (BR) at ( 2, {3 - 3*\hs}) {};  
          
          \path (T2) edge (T1);
          \path (C) edge (T1);
          \path (C) edge (T2);
          
          \path[color=gray] (C) edge (M2);
          \path[color=gray] (M2) edge (M3);
          \path[color=gray] (M3) edge (C);
          
          \path (M2) edge (M1);
          \path (M3) edge (M4);
          
          \path (M1) edge (BL);
          \path (M2) edge (BL);
          \path (BR) edge (M4);
          \path (M3) edge (BR);
        \end{tikzpicture}
         \caption{$C_3\corvf P_2$}
    \end{subfigure}
    \hfill
    \begin{subfigure}[b]{0.3\textwidth}
        \centering
        \begin{tikzpicture}[decoration = {markings,mark=at position .65 with {\arrow{Stealth[length=3mm]}}},
            scale=0.75,
            v/.style={circle,fill,inner sep=1.5pt},
            thick,every edge/.style = {draw, line width=1.pt, postaction=decorate},
        ]
          \pgfmathsetmacro{\hs}{sqrt(3)}  
          
          \node[v] (T1) at (-1,  3) {};            
          \node[v] (T2) at ( 1,  3) {};            
          \node[v] (C)  at ( 0, {3 - \hs}) {};
          
          \node[v] (M1) at (-3, {3 - 2*\hs}) {};   
          \node[v] (M2) at (-1, {3 - 2*\hs}) {};   
          \node[v] (M3) at ( 1, {3 - 2*\hs}) {};   
          \node[v] (M4) at ( 3, {3 - 2*\hs}) {};   
          
          \node[v] (BL) at (-2, {3 - 3*\hs}) {};  
          \node[v] (BR) at ( 2, {3 - 3*\hs}) {};  
          
           \path (T2) edge (T1);
          \path (T1) edge (C);
          \path (T2) edge (C);
          \path[color=gray] (C) edge (M2);
          \path[color=gray] (M2) edge (M3);
          \path[color=gray] (M3) edge (C);
          \path (M1) edge (M2);
          \path (M4) edge (M3);
          
          \path (M1)edge (BL);
          \path (BL) edge (M2);
          \path (BR) edge (M4);
          \path (BR) edge (M3);
        \end{tikzpicture}
         \caption{$C_3 \corvb P_2$}
    \end{subfigure}
    \hfill
    \begin{subfigure}[b]{0.3\textwidth}
        \centering
                \begin{tikzpicture}[decoration = {markings,mark=at position .65 with {\arrow{Stealth[length=3mm]}}},
            scale=0.75,
            v/.style={circle,fill,inner sep=1.5pt},
            thick,every edge/.style = {draw, line width=1.pt, postaction=decorate},
        ]
          \pgfmathsetmacro{\hs}{sqrt(3)}  
          
          \node[v] (T1) at (-1,  3) {};            
          \node[v] (T2) at ( 1,  3) {};            
          \node[v] (C)  at ( 0, {3 - \hs}) {};
          
          \node[v] (M1) at (-3, {3 - 2*\hs}) {};   
          \node[v] (M2) at (-1, {3 - 2*\hs}) {};   
          \node[v] (M3) at ( 1, {3 - 2*\hs}) {};   
          \node[v] (M4) at ( 3, {3 - 2*\hs}) {};   
          
          \node[v] (BL) at (-2, {3 - 3*\hs}) {};  
          \node[v] (BR) at ( 2, {3 - 3*\hs}) {};  
          
        \path (T2) edge (T1);
        \path (T1) edge[bend left=13] (C);
        \path (C) edge[bend left=13] (T1);
        \path (T2) edge[bend left=13] (C);
        \path (C) edge[bend left=13] (T2);

          \path[color=gray] (C) edge (M2);
          \path[color=gray] (M2) edge (M3);
          \path[color=gray] (M3) edge (C);
          \path (M2) edge[bend left=13] (M1);
        \path (M1) edge[bend left=13] (M2);
          \path (M4) edge[bend left=13] (M3);
          \path (M3) edge[bend left=13] (M4);
          
          \path (M1)edge (BL);
          \path (M2) edge[bend left=13] (BL);
          \path (BL) edge[bend left=13] (M2);
          \path (BR) edge (M4);
          \path (M3) edge[bend left=13] (BR);
          \path (BR) edge[bend left=13] (M3);

        \end{tikzpicture}
         \caption{$C_3\corv P_2$}
    \end{subfigure}
    \caption{The forward-vertex-corona, backward-vertex-corona and symmetric-vertex-corona of the directed cycle $C_3$ (indicated by gray arcs) and directed path $P_2$.}
    \label{dual1.1}
\end{figure}
\cref{dual1.1} illustrates the various vertex-coronas for digraphs when $D_1$ is the directed cycle $C_3$ and $D_2$ is the directed path $P_2$.
The symmetric-vertex-corona coincides with the classical corona of two graphs when $D_1$ and $D_2$ are symmetric digraphs (i.e., undirected graphs).
We are not aware of any resources that attempt to generalize the corona product to digraphs with the following exception: In \cite{hasyyati2021characteristic}, the characteristic polynomial of the anti-adjacency matrix for the digraph $C_n\corvb\overline{K_r}$ is analyzed, where $C_n$ is the directed $n$-cycle and $\overline{K_r}$ is the empty digraph with $r$ vertices and no arcs; the authors refer to this digraph as a ``directed unicyclic corona graph''.

Given any digraphs $D_1$ and $D_2$, the digraphs $D_1 \corvf D_2$ and $D_1 \corvb D_2$ are not strongly connected.
In such cases, the adjacency matrix is block triangular (up to permutation) and its spectrum is the union of the spectra of the diagonal blocks.  
For the Laplacian and signless Laplacian, additional arcs between strong components affect $D_{\textnormal{out}}(D)$, so extra care is required.  
We therefore focus on corona products that may yield strongly connected digraphs and thus restrict our attention to the symmetric-vertex-corona in this section.  
We begin with a remark on when the symmetric-vertex-corona is strongly connected.

\begin{remark}\rm
Let $D_1$ and $D_2$ be digraphs. 
Then $D_1\overleftrightarrow{\circ} D_2$ is strongly connected if and only if $D_1$ is strongly connected. 
\end{remark}

To describe the adjacency (resp. Laplacian and signless Laplacian) matrix of the symmetric-vertex-corona of two digraphs $D_1$ and $D_2$, we adopt the vertex labelling described in \cite[Section~2]{Mcleman}.
In particular, 
first choose an arbitrary ordering $v_1,\ldots,v_{n_2}$ of the vertices of $D_2$.
We then label the vertices in the copy of $D_1$ by $1,2,\ldots,n_1$, and for $1\leq i\leq n_1$ and $1\leq k\leq n_2$, label the vertex in the $i$th copy of $D_2$ corresponding to $v_k$ by the label $i+n_1k$.
Under this labelling, we have
\[
A(D_1\corv D_2)=
\left[\begin{array}{c|c}
A(D_1) & \textbf{1}_{n_2}^T\otimes I_{n_1} \\
\hline
\textbf{1}_{n_2}\otimes I_{n_1} & A(D_2)\otimes I_{n_1}\\
\end{array}\right],~~
L(D_1\corv D_2)=
\left[\begin{array}{c|c}
L(D_1)+n_2I_{n_1} & -\textbf{1}_{n_2}^T\otimes I_{n_1} \\
\hline
-\textbf{1}_{n_2}\otimes I_{n_1} & (L(D_2)+I_{n_2})\otimes I_{n_1}\\
\end{array}\right], 
\] 
\[
\text{and}~~
Q(D_1\corv D_2)=
\left[\begin{array}{c|c}
Q(D_1)+n_2I_{n_1} & \textbf{1}_{n_2}^T\otimes I_{n_1} \\
\hline
\textbf{1}_{n_2}\otimes I_{n_1} & (Q(D_2)+I_{n_2})\otimes I_{n_1}\\
\end{array}\right].
\] 
In \cite[Theorem 2]{Mcleman}, the adjacency characteristic polynomial of the corona of two graphs is computed using the coronal. 
This result is extended to the Laplacian in \cite{LIUu} (note that the $L$-coronal is defined slightly differently in \cite{LIUu}) and to the signless Laplacian in \cite{cui2012spectrum}. 
Since the proofs do not rely on symmetry, they apply to digraphs and to general matrices. 
The core argument yields the general relation below, which follows by the same reasoning as \cite[Theorem~2]{Mcleman} using properties of Kronecker products and \cref{schur}. 
We leave the proof of this to the reader.

\begin{lemma}\label[lemma]{lm:kronschur}
Suppose
\[
M=\left[\begin{array}{cc}
M_1 & \pm\textbf{1}_{n_2}^T\otimes B_{1} \\
\pm\textbf{1}_{n_2}\otimes B_2 & M_2\otimes I_{r}\\\end{array}\right],
\]
where $M_1$, $M_2$, $B_1$ and $B_2$ are $n_1\times n_1$, $n_2\times n_2$, $n_1\times r$ and $r\times n_1$ matrices, respectively.
Then
\[
\begin{aligned}
f_M(\lambda) 
&= \big[f_{M_2}(\lambda)\big]^{r} \det\left(\lambda I_{n_1}-M_1- \chi_{M_2}(\lambda)\,B_1B_2\right).
\end{aligned}
\]
Furthermore, in the case that $B_1=B_2=I_{n_1}$, we have
$
\begin{aligned}
f_M(\lambda) 
&= \big[f_{M_2}(\lambda)\big]^{n_1} f_{M_1}\left(\lambda - \chi_{M_2}(\lambda)\right).
\end{aligned}
$
\end{lemma}

\cref{cor_addI}(i) and \cref{lm:kronschur} give the main result of this section.

\begin{theorem}\label[theorem]{coronal-thm}
Let $D_1$ and $D_2$ be digraphs with $n_1$ and $n_2$ vertices, respectively.
Then
\begin{enumerate}[(i)]
    \item $\displaystyle f_{A(D_1\corv D_2)}(\lambda) = \big[f_{A(D_2)}(\lambda)\big]^{n_1} f_{A(D_1)}\left(\lambda-\chi_{A(D_2)}(\lambda)\right)$, 
    \item $\displaystyle f_{L(D_1\corv D_2)}(\lambda) =\big[f_{L(D_2)}(\lambda-1)\big]^{n_1} f_{L(D_1)}\left(\frac{\lambda^2-(n_2+1)\lambda}{\lambda-1}\right)$, 
    \item $\displaystyle f_{Q(D_1\corv D_2)}(\lambda) =\big[f_{Q(D_2)}(\lambda-1)\big]^{n_1} f_{Q(D_1)}\left(\lambda-n_2-\chi_{Q(D_2)}(\lambda-1)\right)$.
\end{enumerate}
In particular, the adjacency (resp. signless Laplacian) spectrum of $D_1\corv D_2 $ is completely determined by the adjacency (resp.  signless Laplacian) characteristic polynomials of $D_1$ and $D_2$ and the $A$-coronal (resp. $Q$-coronal) of $D_2$, whereas,
the Laplacian spectrum of $D_1\corv D_2 $ is completely determined by the Laplacian characteristic polynomials of $D_1$ and $D_2$.
\end{theorem}

\begin{proof}
For the first equation, the result follows directly from \cref{lm:kronschur} by taking $M_1=A(D_1)$, $M_2=A(D_2)$ and $B_1=B_2=I_{n_1}$.
For the second equation, take
$M_1=L(D_1)+n_2I_{n_1}$, $M_2=L(D_2)+I_{n_2}$ and $B_1=B_2=I_{n_1}$ in \cref{lm:kronschur}, then apply \cref{cor_addI} and \cref{conoal_cal} to obtain
\begin{align*}  
f_{L(D_1\corv D_2)}(\lambda)
&= \big[f_{L(D_2)+I_{n_2}}(\lambda)\big]^{n_1} f_{L(D_1)+n_2I_{n_1}}\left(\lambda - \chi_{L(D_2)+I_{n_2}}(\lambda)\right)\\
&=\big[f_{L(D_2)}(\lambda-1)\big]^{n_1} f_{L(D_1)}\left(\lambda-\chi_{L(D_2)+I_{n_2}}(\lambda)-n_2\right)\\
&=\big[f_{L(D_2)}(\lambda-1)\big]^{n_1} f_{L(D_1)}\left(\lambda-\chi_{L(D_2)}(\lambda-1)-n_2\right)\\
&=\big[f_{L(D_2)}(\lambda-1)\big]^{n_1} f_{L(D_1)}\left(\lambda-\frac{n_2}{\lambda-1}-n_2\right).
\end{align*}  
The third equation is derived by an identical argument: take $M_1=Q(D_1)+n_2I_{n_1}$, 
$M_2=Q(D_2)+I_{n_2}$ and $B_1=B_2=I_{n_1}$ in \cref{lm:kronschur} and apply \cref{cor_addI}.
\end{proof}

In \cite{barik2007spectrum} and \cite[Proposition 6]{Mcleman}, the adjacency spectrum of the corona of a graph and a regular graph is determined and extended to the Laplacian in \cite{LIUu} and to the signless Laplacian in \cite{cui2012spectrum}.
These results also hold for out-regular digraphs. 

\begin{corollary}{\rm}\label[corollary]{r-regular-coronal}
Let $D_1$ and $D_2$ be digraphs with $n_1$ and $n_2$ vertices, respectively, and suppose that $D_2$ is $r$-out-regular and strongly connected.
Then the following statements hold.
\begin{enumerate}[(i)] 
\item The spectrum of $A(D_1 \corv D_2)$ consists of all eigenvalues of $A(D_2)$ not equal to $r$, each with multiplicity $n_1$, together with two multiplicity-one eigenvalues of the form $$\frac{1}{2}\left(\mu + r \pm \sqrt{(r - \mu)^2 + 4n_2}\right)$$ for each eigenvalue $\mu$ of $A(D_1)$.

\item The spectrum of $L(D_1 \corv D_2)$ consists of all values $\delta+1$, for each eigenvalue $\delta$ of $L(D_2)$ not equal to $0$, each with multiplicity $n_1$, along with two multiplicity-one eigenvalues of the form $$\frac{1}{2}\left(\mu + n_2 + 1 \pm \sqrt{(\mu + n_2 + 1)^2 - 4\mu}\right)$$ for each eigenvalue $\mu$ of $L(D_1)$. 

\item The spectrum of $Q(D_1 \corv D_2)$ consists of all values $\delta+1$, for each eigenvalue $\delta$ of $Q(D_2)$ not equal to $2r$, each with multiplicity $n_1$, along with two multiplicity-one eigenvalues of the form $$\frac{1}{2}\left(\mu+n_2+2r+1\pm\sqrt{((\mu+n_2)-(2r+1))^2+4n_2}\right)$$
for each eigenvalue $\mu$ of $Q(D_1)$.  
\end{enumerate}
\end{corollary} 
 
\begin{proof}
Since $D_2$ is $r$-out-regular and strongly connected, $\lambda = r$ is a simple root of $f_{A(D_2)}(\lambda)$, and by \cref{conoal_cal}, $\chi_{A(D_2)}(\lambda) = n_2 / (\lambda - r)$ has a single pole at $\lambda = r$. 
The result now follows from \cref{coronal-thm} by solving $\lambda - \chi_{A(D_2)}(\lambda) = \mu$ for each eigenvalue $\mu$ of $A(D_1)$.
The same reasoning applies to $L(D_2)$ and $Q(D_2)$. In these cases, the coronals are given by $\chi_{L(D_2)}(\lambda - 1) = n_2 / (\lambda - 1)$ and $\chi_{Q(D_2)}(\lambda - 1) = n_2 / (\lambda - 1 - 2r)$, each with a single pole at $\lambda = 1$ and $\lambda = 2r + 1$, respectively. These correspond to simple roots of $f_{L(D_2)}(\lambda - 1)$ and $f_{Q(D_2)}(\lambda - 1)$ since $D_2$ is strongly connected (e.g., see \cite[Proposition 4.5]{veerman2020primer} that shows $0$ is a simple eigenvalue of $L(D_2)$). The result follows from solving
$(\lambda^2-(n_2+1)\lambda)/(\lambda-1) = \mu$
and $\lambda - n_2 - \chi_{Q(D_2)}(\lambda - 1) = \mu$
for each eigenvalue $\mu$ of $L(D_1)$ and $Q(D_1)$, respectively.
This gives
\[
\lambda=\frac{1}{2}\left(\mu+n_2+2r+1\pm\sqrt{((\mu+n_2)+(2r+1))^2-4(\mu(2r+1)+2rn_2)}\right),
\]
which simplifies to the stated equation.
\end{proof}


When the coronal of the digraph $D_2$ has a simple form, \cref{coronal-thm} can be used to obtain simple formulas for the characteristic polynomial of the symmetric-vertex-corona of $D_1$ and $D_2$.
For example, 
if $D_2=P_{n_2}$ is the directed path on $n_2\geq 1$ vertices,
then \cref{prop:coronoal_path} gives
\[
f_{A(D_1\corv P_{n_2})}(\lambda)
=\lambda^{n_1n_2}f_{A(D_1)}\left(\lambda-\frac{{n_2}\lambda^{{n_2}+1}-(n_2+1)\lambda^{n_2}+1}{\lambda^{n_2}(\lambda-1)^2}\right).
\]
In some cases, the spectrum can also be derived using the same approach as in \cref{r-regular-coronal}.
This also applies to the families of digraphs described in Section~\ref{sec:cor_comp}, but we omit these computations in this paper.



Some of the results in this section can be extended to real matrices and have connections to structured matrices and systems of linear second-order ordinary differential equations (see \cite{berliner2023inertias} and references therein).
To be specific, when $D_2$ is a single vertex (i.e., $D_2=P_1$), then (up to permutation similarity)
\[
A(D_1\corv P_1)=
\left[\begin{array}{c|c}
A(D_1) & I_{n_1} \\
\hline
I_{n_1} & O_{n_1}\\
\end{array}\right].
\]
The spectrum and inertia of matrices of this form has been studied in the context of sign patterns where $A(D_1)$ is replaced by a real matrix whose entries have specified signs and the $(1,2)$-block of $A(D_1\corv P_1)$ is a positive diagonal matrix.
In this context, \cite[Lemma~2.1]{berliner2023inertias} gives a generalization to \cref{r-regular-coronal}(i) in the case when $D_2=P_1$.

\section{Arc Corona}\label{sec:arccorona}
Hou and Shiu \cite{Hou} first defined the edge-corona of two graphs and analyzed the spectra of their adjacency and Laplacian matrices.
This was extended to include the signless Laplacian spectrum in \cite{wang2013signless} and to weighted networks in \cite{liu2020generalized}.
Specifically, for graphs $G_1$ and $G_2$ with $G_1$ having $m_1$ edges, the \emph{edge-corona} $G_1 \diamond G_2$ is formed by taking one copy of $G_1$ and $m_1$ copies of $G_2$, then joining both endpoints of each edge $e_k = uv$ of $G_1$ to every vertex in the $k$th copy of $G_2$. 
We propose three extensions of this concept to digraphs that potentially give ``arc-corona'' products that are strongly connected.

\begin{defn}\rm\label{arc-corona} 
Let $D_1$ and $D_2$ be digraphs with $D_1$ having $m_1$ arcs and its underlying graph $U(D_1)$ having $m_1'$ edges. 
The \defin{forward-arc-corona} $D_1 \arcf D_2$ is formed by taking one copy of $D_1$ and $m_1$ copies of $D_2$, and for each arc $a_k = uv$ in $E(D_1)$ and each vertex $w$ in the $k$th copy of $D_2$, adding arcs $uw$ and $wv$. 
The \defin{backward-arc-corona} $D_1 \arcb D_2$ instead adds arcs $vw$ and $wu$.
The \defin{symmetric-arc-corona} $D_1 \arcs D_2$ is formed by taking one copy of $D_1$ and $m_1'$ copies of $D_2$, and for each corresponding edge $a_k=uv$ in $E(U(D_1))$ and each vertex $w$ in the $k$th copy of $D_2$, adding arcs $uw$, $wu$, $vw$ and $wv$.
\end{defn}

See \cref{fig:arccorona_connections} for an illustration of the arc configurations described in \cref{arc-corona} and see \cref{{dual1.2}} illustrating $C_3 \arcf P_2$, $C_3 \arcb P_2$ and $C_3 \arcs P_2$.
We would like to emphasize that for the symmetric-arc-corona, each symmetric edge in $E(D_1)$ (i.e., pair of arcs $uv$ and $vu$) adds one copy of $D_2$ to $D_1 \arcs D_2$ rather than one copy for each arc.
There are two reasons for this.
First, when $D_1$ and $D_2$ are symmetric digraphs, this definition coincides with the edge-corona of two graphs, and second, in \cref{thm:arccor_char} we observe that the signless Laplacian of $U(D_1)$ naturally appears in the characteristic polynomial computation for $A(D_1\arcs D_2)$ under this convention.
\begin{figure}[!ht]
    \centering
    \begin{subfigure}[b]{0.24\textwidth}
        \centering
        \begin{tikzpicture}[decoration = {markings,mark=at position .65 with {\arrow{Stealth[length=3mm]}}},
            dot/.style = {circle, fill, inner sep=1.3pt, node contents={},label=#1},
            every edge/.style = {draw, line width=1.pt, postaction=decorate},
            scale=0.75]
            \node (u) at (0,0) [dot=left:$u$];
            \node (v) at (0,2) [dot=left:$v$];
            \node (w) at (1.414213,1) [dot=right:$w$];
            \path (u) edge (v);
            \path (u) edge (w);
            \path (w) edge (v);
        \end{tikzpicture}
        \caption{$P_2 \arcf P_1$}
    \end{subfigure}
    \hfill
    \begin{subfigure}[b]{0.24\textwidth}
        \centering
        \begin{tikzpicture}[decoration = {markings,mark=at position .65 with {\arrow{Stealth[length=3mm]}}},
            dot/.style = {circle, fill, inner sep=1.3pt, node contents={},label=#1},
            every edge/.style = {draw, line width=1.pt, postaction=decorate},
            scale=0.75]
            \node (u) at (0,0) [dot=left:$u$];
            \node (v) at (0,2) [dot=left:$v$];
            \node (w) at (1.414213,1) [dot=right:$w$];
            \path (u) edge (v);
            \path (w) edge (u);
            \path (v) edge (w);
        \end{tikzpicture}
        \caption{$P_2 \arcb P_1$}
    \end{subfigure}
    \hfill
    \begin{subfigure}[b]{0.24\textwidth}
        \centering
        \begin{tikzpicture}[decoration = {markings,mark=at position .65 with {\arrow{Stealth[length=3mm]}}},
            dot/.style = {circle, fill, inner sep=1.3pt, node contents={},label=#1},
            every edge/.style = {draw, line width=1.pt, postaction=decorate},
            scale=0.75]
            \node (u) at (0,0) [dot=left:$u$];
            \node (v) at (0,2) [dot=left:$v$];
            \node (w) at (1.414213,1) [dot=right:$w$];
            \path (u) edge (v);
            \path (u) edge[bend left=13] (w);
            \path (w) edge[bend left=13] (u);
            \path (v) edge[bend left=13] (w);
            \path (w) edge[bend left=13] (v);
        \end{tikzpicture}
        \caption{$P_2 \arcs P_1$}
    \end{subfigure}
    \begin{subfigure}[b]{0.24\textwidth}
        \centering
        \begin{tikzpicture}[decoration = {markings,mark=at position .65 with {\arrow{Stealth[length=3mm]}}},
            dot/.style = {circle, fill, inner sep=1.3pt, node contents={},label=#1},
            every edge/.style = {draw, line width=1.pt, postaction=decorate},
            scale=0.75]
            \node (u) at (0,0) [dot=left:$u$];
            \node (v) at (0,2) [dot=left:$v$];
            \node (w) at (1.414213,1) [dot=right:$w$];
            \path (u) edge[bend left=13] (v);
            \path (v) edge[bend left=13] (u);
            \path (u) edge[bend left=13] (w);
            \path (w) edge[bend left=13] (u);
            \path (v) edge[bend left=13] (w);
            \path (w) edge[bend left=13] (v);
        \end{tikzpicture}
        \caption{$C_2 \arcs P_1$}
    \end{subfigure}    
    \caption{Possible connections between an arc $uv$ in $D_1=P_2$ and a vertex $w$ in $D_2=P_1$ in forming the arc corona. In (d), the symmetric edge $uv$ adds one copy of $D_2$ to the symmetric-arc-corona.}
    \label{fig:arccorona_connections}
\end{figure}
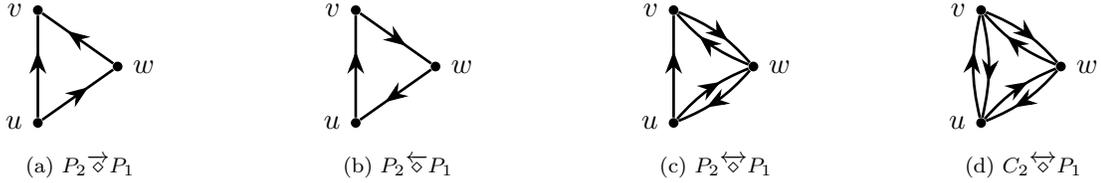


\begin{figure}[!ht]
    \centering  
    \begin{subfigure}[b]{0.3\textwidth}
        \centering
                \begin{tikzpicture}[decoration = {markings,mark=at position .72 with {\arrow{Stealth[length=3mm]}}},
            scale=0.75,
            v/.style={circle,fill,inner sep=1.5pt},
            thick,every edge/.style = {draw, line width=1.pt, postaction=decorate},
        ]
          \pgfmathsetmacro{\hs}{sqrt(3)}  
          \node[v] (V1)  at ( 0, {3 - \hs}) {};
          \node[v] (V2) at ( 1, {3 - 2*\hs}) {};
          \node[v] (V3) at (-1, {3 - 2*\hs}) {};   
          \node[v] (W1) at ({-\hs}, {4-\hs}) {}; 
          \node[v] (W2) at ({-1-\hs}, {4-2*\hs}) {};
          \node[v] (W3) at (-1, {1 - 2*\hs}) {}; 
          \node[v] (W4) at ( 1, {1 - 2*\hs}) {}; 
          \node[v] (W6) at (\hs, 4-\hs) {};            
          \node[v] (W5) at ({1+\hs}, {4-2*\hs}) {}; 
          
          \path[color=gray] (V1) edge (V3);
          \path[color=gray] (V3) edge (V2);
          \path[color=gray] (V2) edge (V1);          
          \path (W1) edge (W2);
          \path (W3) edge (W4);
          \path (W5) edge (W6);
          \path (V1) edge (W1);
          \path (W1) edge (V3);
          \path (V1) edge (W2);
          \path (W2) edge (V3);
          \path (V3) edge (W3);
          \path (W3) edge (V2);
          \path (V3) edge (W4);
          \path (W4) edge (V2);          
          \path (V2) edge (W5);
          \path (W5) edge (V1);
          \path (V2) edge (W6);
          \path (W6) edge (V1);          
        \end{tikzpicture}
         \caption{$C_3\corvf P_2$}
    \end{subfigure}
    \hfill
    \begin{subfigure}[b]{0.3\textwidth}
        \centering
                \begin{tikzpicture}[decoration = {markings,mark=at position .72 with {\arrow{Stealth[length=3mm]}}},
            scale=0.75,
            v/.style={circle,fill,inner sep=1.5pt},
            thick,every edge/.style = {draw, line width=1.pt, postaction=decorate},
        ]
          \pgfmathsetmacro{\hs}{sqrt(3)}  
          \node[v] (V1)  at ( 0, {3 - \hs}) {};
          \node[v] (V2) at ( 1, {3 - 2*\hs}) {};
          \node[v] (V3) at (-1, {3 - 2*\hs}) {};   
          \node[v] (W1) at ({-\hs}, {4-\hs}) {}; 
          \node[v] (W2) at ({-1-\hs}, {4-2*\hs}) {};
          \node[v] (W3) at (-1, {1 - 2*\hs}) {}; 
          \node[v] (W4) at ( 1, {1 - 2*\hs}) {}; 
          \node[v] (W6) at (\hs, 4-\hs) {};            
          \node[v] (W5) at ({1+\hs}, {4-2*\hs}) {}; 
          
          \path[color=gray] (V1) edge (V3);
          \path[color=gray] (V3) edge (V2);
          \path[color=gray] (V2) edge (V1);          
          \path (W1) edge (W2);
          \path (W3) edge (W4);
          \path (W5) edge (W6);
          \path (W1) edge (V1);
          \path (V3) edge (W1);
          \path (W2) edge (V1);
          \path (V3) edge (W2);
          \path (W3) edge (V3);
          \path (V2) edge (W3);
          \path (W4) edge (V3);
          \path (V2) edge (W4);          
          \path (W5) edge (V2);
          \path (V1) edge (W5);
          \path (W6) edge (V2);
          \path (V1) edge (W6);          
        \end{tikzpicture}
         \caption{$C_3 \corvb P_2$}
    \end{subfigure}
    \hfill
    \begin{subfigure}[b]{0.3\textwidth}
        \centering
                \begin{tikzpicture}[decoration = {markings,mark=at position .75 with {\arrow{Stealth[length=3mm]}}},
            scale=0.75,
            v/.style={circle,fill,inner sep=1.5pt},
            thick,every edge/.style = {draw, line width=1.pt, postaction=decorate},
        ]
          \pgfmathsetmacro{\hs}{sqrt(3)}
          \pgfmathsetmacro{\ebend}{10}
          \node[v] (V1)  at ( 0, {3 - \hs}) {};
          \node[v] (V2) at ( 1, {3 - 2*\hs}) {};
          \node[v] (V3) at (-1, {3 - 2*\hs}) {};   
          \node[v] (W1) at ({-\hs}, {4-\hs}) {}; 
          \node[v] (W2) at ({-1-\hs}, {4-2*\hs}) {};
          \node[v] (W3) at (-1, {1 - 2*\hs}) {}; 
          \node[v] (W4) at ( 1, {1 - 2*\hs}) {}; 
          \node[v] (W6) at (\hs, 4-\hs) {};            
          \node[v] (W5) at ({1+\hs}, {4-2*\hs}) {}; 
          
          \path[color=gray] (V1) edge (V3);
          \path[color=gray] (V3) edge (V2);
          \path[color=gray] (V2) edge (V1);          
          \path (W1) edge (W2);
          \path (W3) edge (W4);
          \path (W5) edge (W6);
          \path (V1) edge[bend left=\ebend] (W1);
          \path (W1) edge[bend left=\ebend] (V3);
          \path (V1) edge[bend left=\ebend] (W2);
          \path (W2) edge[bend left=\ebend] (V3);
          \path (V3) edge[bend left=\ebend] (W3);
          \path (W3) edge[bend left=\ebend] (V2);
          \path (V3) edge[bend left=\ebend] (W4);
          \path (W4) edge[bend left=\ebend] (V2);          
          \path (V2) edge[bend left=\ebend] (W5);
          \path (W5) edge[bend left=\ebend] (V1);
          \path (V2) edge[bend left=\ebend] (W6);
          \path (W6) edge[bend left=\ebend] (V1);            
          \path (W1) edge[bend left=\ebend] (V1);
          \path (V3) edge[bend left=\ebend] (W1);
          \path (W2) edge[bend left=\ebend] (V1);
          \path (V3) edge[bend left=\ebend] (W2);
          \path (W3) edge[bend left=\ebend] (V3);
          \path (V2) edge[bend left=\ebend] (W3);
          \path (W4) edge[bend left=\ebend] (V3);
          \path (V2) edge[bend left=\ebend] (W4);          
          \path (W5) edge[bend left=\ebend] (V2);
          \path (V1) edge[bend left=\ebend] (W5);
          \path (W6) edge[bend left=\ebend] (V2);
          \path (V1) edge[bend left=\ebend] (W6);          
        \end{tikzpicture}
         \caption{$C_3\corv P_2$}
    \end{subfigure}
    \caption{The forward-vertex-corona, backward-vertex-corona and symmetric-vertex-corona of the directed cycle $C_3$ (indicated with gray arcs) and directed path $P_2$.}
    \label{dual1.2}
\end{figure}
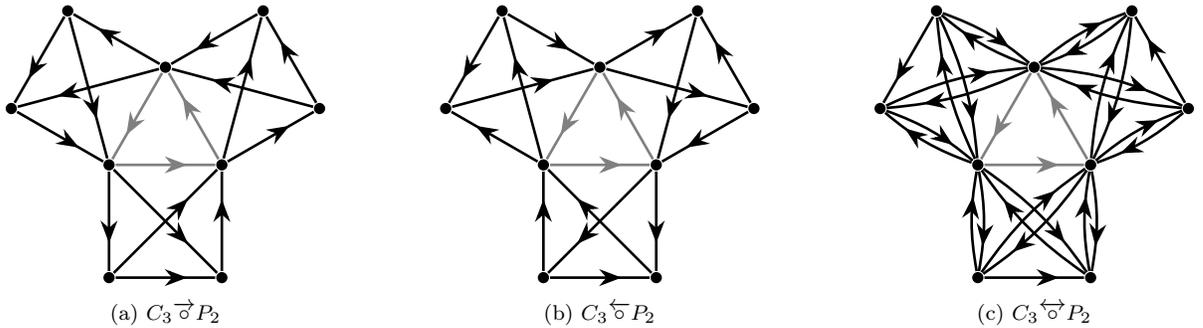

Alternate definitions that necessarily produce digraphs that are not strongly connected are not considered here. 
Furthermore, we note that $D_1 \arcf D_2$ (resp. $D_1 \arcb D_2$ and $D_1 \arcs D_2$) are well-defined in the sense that different vertex orderings of $D_1$ and $D_2$ produce arc-corona adjacency matrices that are permutation similar.
In general, the digraphs $D_1 \arcf D_2$ and $D_1 \arcb D_2$ are nonisomorphic (e.g., see \cref{fig:arccorona_connections} comparing $P_2 \arcf P_1$ and $P_2 \arcb P_1$),
however, if $D_1$ is a symmetric digraph, then $D_1\arcf D_2\cong D_1\arcb D_2$ holds for any digraph $D_2$.

The following remark establishes conditions under which the arc-corona of two digraphs produces a strongly connected digraph.
It follows by observing that for an arc $uv$ in $D_1$, the construction of $D_1 \arcf D_2$ does not introduce any new paths from $v$ to $u$, whereas $D_1 \arcb D_2$ and $D_1 \arcs D_2$ explicitly add such a path.

\begin{remark}\label[remark]{str}\rm
Let $D_1$ and $D_2$ be digraphs.
Then:
\begin{enumerate}[(i)]
    \item $D_1\arcf D_2$ is strongly connected if and only if $D_1$ is strongly connected.
    \item $D_1\arcb D_2$ is strongly connected if and only if the underlying graph of $D_1$ is connected.
    \item $D_1\arcs D_2$ is strongly connected if and only if the underlying graph of $D_1$ is connected.
\end{enumerate}
\end{remark}

\subsection{Adjacency spectrum of arc coronas}
We now give formulas for the adjacency characteristic polynomials for the three arc-corona constructions.

\begin{theorem}\label[theorem]{thm:arccor_char}
Let $ D_1 $ and $ D_2 $ be two digraphs on $ n_1 $ and $ n_2 $ vertices, respectively.
Suppose $D_1$ has $m_1$ arcs and its underlying digraph 
$G_1=U(D_1)$ has $m_1'$ edges. 
Then 
\begin{enumerate}[(i)]
    \item $\displaystyle  f_{A(D_1 \arcf D_2)} (\lambda)=\big[f_{A(D_2)}(\lambda)\big]^{m_1}\left(1+\chi_{A(D_2)}(\lambda)\right)^{n_1}f_{A(D_1)}\left(\frac{\lambda}{1+\chi_{A(D_2)}(\lambda)}\right)$, 
    \item $\displaystyle f_{A(D_1 \arcb D_2)}(\lambda)= \big[f_{A(D_2)}(\lambda)\big]^{m_1} \det\left(\lambda I_{n_1}-A(D_1)-\chi_{A(D_2)}(\lambda)\,A(\overleftarrow{D}_1)\right)$, 
    \item $\displaystyle f_{A(D_1 \arcs D_2)}(\lambda)= \big[f_{A(D_2)}(\lambda)\big]^{m_1'} \det\left(\lambda I_{n_1}-A(D_1)-\chi_{A(D_2)}(\lambda)\,Q(G_1)\right)$.
\end{enumerate}
In particular, the adjacency spectrum of $D_1 \arcf D_2$ is completely determined by the adjacency characteristic polynomials of $D_1$ and $D_2$, and the $A$-coronal of $D_2$.
\end{theorem}
 
\begin{proof}
The matrices $A(D_1\arcf D_2)$, $A(D_1\arcb D_2)$ and $A(D_1\arcs D_2)$ are permutation similar to
\[
\left[\begin{array}{cc}
A(D_1) & \textbf{1}_{n_2}^T\otimes B_{\rm out}(D_1) \\
\textbf{1}_{n_2}\otimes B_{\rm in}(D_1)^T & A(D_2)\otimes I_{m_1}\\
\end{array}\right],~
\left[\begin{array}{cc}
A(D_1) & \textbf{1}_{n_2}^T\otimes B_{\rm in}(D_1) \\
\textbf{1}_{n_2}\otimes B_{\rm out}(D_1)^T & A(D_2)\otimes I_{m_1}\\
\end{array}\right],
\]
\[\text{and}~
\left[\begin{array}{cc}
A(D_1) & \textbf{1}_{n_2}^T\otimes B(G_1) \\
\textbf{1}_{n_2}\otimes B(G_1)^T & A(D_2)\otimes I_{m_1'}\\
\end{array}\right],
\] 
respectively,
where $B_{\rm in}(D_1)$ and $B_{\rm out}(D_1)$ are the $n_1\times m_1$ in-incidence and out-incidence matrices of $D_1$, and $B(G_1)$ is the $n_1\times m_1'$ incidence matrix of the graph $G_1$.
We apply \cref{lm:kronschur} to each of these matrices and note that
\cref{bb} gives $B_{\rm out}(D_1)B_{\rm in}(D_1)^T=A(D_1)$ and 
$B_{\rm in}(D_1)B_{\rm out}(D_1)^T=(B_{\rm out}(D_1)B_{\rm in}(D_1)^T)^T=A(D_1)^T=A(\overleftarrow{D}_1)$, and for the graph $G_1$ we have $B(G_1)B(G_1)^T=Q(G_1)$.
The second and third equations immediately follow, whereas for the first equation, we further note that
\begin{align*}
\det\left(\lambda I_{n_1}-A(D_1)-\chi_{A(D_2)}(\lambda)\,A(D_1)\right)
&=
\det\left(\lambda I_{n_1}-(1+\chi_{A(D_2)}(\lambda))A(D_1)\right)\\
&=
\left(1+\chi_{A(D_2)}(\lambda)\right)^{n_1}\det\left(\frac{\lambda}{1+\chi_{A(D_2)}(\lambda)} I_{n_1}-A(D_1)\right).\\
&=\left(1+\chi_{A(D_2)}(\lambda)\right)^{n_1}
f_{A(D_1)}\left(\frac{\lambda}{1+\chi_{A(D_2)}(\lambda)}\right),
\end{align*}
as required.
\end{proof}

When $D_1$ and $D_2$ are symmetric digraphs (graphs), we have the following consequence.
\begin{corollary}
Let $G_1$ and $G_2$ be two graphs on $n_1$ and $n_2$ vertices, and $m_1$ and $m_2$ edges, respectively. 
Suppose $Q(G_1)$ is the signless Laplacian matrix of $G_1$. 
Then 
\[
f_{A(G_1 \diamond G_2)}(\lambda)= \big[f_{A(G_2)}(\lambda)\big]^{m_1} \det\left(\lambda I_{n_1}-A(G_1)-\chi_{A(G_2)}(\lambda)\,Q(G_1)\right).
\]
\end{corollary}

Our interest is in cases where the adjacency spectrum of the arc-corona of digraphs $D_1$ and $D_2$ can be determined from the adjacency characteristic polynomials of $D_1$ and $D_2$ together with the $A$-coronal of $D_2$ (and possibly the $A$-coronal of $D_1$).
\cref{thm:arccor_char} verifies this is true for the forward-arc-corona for arbitrary digraphs $D_1$ and $D_2$. 
As mentioned previously, in the case that $D_1$ is a symmetric digraph, we have $D_1\arcf D_2\cong D_1\arcb D_2$ for any digraph $D_2$.
This gives the following corollary for the backward-arc-corona.

\begin{corollary}\label[corollary]{rem:arcb_graph}
Let $D_1$ and $D_2$ be two digraphs on $n_1$ and $n_2$ vertices, and $m_1$ and $m_2$ arcs, respectively. 
If $D_1$ is a symmetric digraph, then 
\[
f_{A(D_1 \arcb D_2)} (\lambda) =\big[f_{A(D_2)}(\lambda)\big]^{m_1}
\left(1+\chi_{A(D_2)}(\lambda)\right)^{n_1}
f_{A(D_1)}\left(\frac{\lambda}{1+\chi_{A(D_2)}(\lambda)}\right).
\]
In particular, when $D_1$ is a symmetric digraph, the adjacency spectrum of $D_1 \arcb D_2$ is completely determined by the adjacency characteristic polynomials of $D_1$ and $D_2$, and the $A$-coronal of $D_2$.
\end{corollary}

Another class of digraphs to consider for the backward-arc-corona are tournaments: a \emph{tournament} is a digraph $D$ with $n$ vertices where $A(D)+A(\overleftarrow{D})=J_n-I_n$ holds.

\begin{corollary}
Let $ D_1 $ and $ D_2 $ be two digraphs on $ n_1 $ and $ n_2 $ vertices, and $ m_1 $ and $ m_2 $ arcs, respectively. 
Suppose $D_1$ is a tournament and let $\chi_1(\lambda)=\chi_{A(D_1)}(\lambda)$ and $\chi_2(\lambda)=\chi_{A(D_2)}(\lambda)$.
Then $f_{A(D_1 \arcb D_2)} (\lambda)$ is 
\begin{align*}  
\big[f_{A(D_2)}(\lambda)\big]^{m_1}
\big[1-\chi_2(\lambda)\big]^{n_1-1}
f_{A(D_1)}\left(\frac{\lambda+\chi_2(\lambda)}{1-\chi_2(\lambda)}\right)
\left(1-\chi_2(\lambda)-\chi_2(\lambda)\,\chi_{1}\left(\frac{\lambda+\chi_2(\lambda)}{1-\chi_2(\lambda)}\right)\right).
\end{align*}
In particular, when $D_1$ is a tournament, the adjacency spectrum of $D_1 \arcb D_2$ is completely determined by the adjacency characteristic polynomials of $D_1$ and $D_2$, and the $A$-coronals of $D_1$ and $D_2$.
\end{corollary}

\begin{proof}
Since $D_1$ is a tournament, using properties of the determinant and \cref{sylv}(i) we obtain

$\det\left(\lambda I_{n_1}-A(D_1)-\chi_{2}(\lambda)\,A(\overleftarrow{D}_1)\right)$\vspace{-0.75em}
\begin{align*}
\hspace{2.5em}&=
\det\big(\lambda I_{n_1}-A(D_1)-\chi_2(\lambda)\,(J_{n_1}-I_{n_1}-A(D_1))\big)\\
&=
\det\big((\lambda+\chi_2(\lambda)) I_{n_1}-(1-\chi_2(\lambda))A(D_1)-\chi_2(\lambda)\,J_{n_1}\big)\\
&=\big[1-\chi_2(\lambda)\big]^{n_1}
\det\left(\frac{\lambda+\chi_2(\lambda)}{1-\chi_2(\lambda)} I_{n_1}-A(D_1)-\frac{\chi_2(\lambda)}{1-\chi_2(\lambda)}\,J_{n_1}\right)\\
&=\big[1-\chi_2(\lambda)\big]^{n_1}
\det\left(\frac{\lambda+\chi_2(\lambda)}{1-\chi_2(\lambda)} I_{n_1}-A(D_1)\right)
\left(1-\frac{\chi_2(\lambda)}{1-\chi_2(\lambda)}\,\textbf{1}_n^T\left(\frac{\lambda+\chi_2(\lambda)}{1-\chi_2(\lambda)} I_{n_1}-A(D_1)\right)^{-1}\textbf{1}_n\right)\\
&=\big[1-\chi_2(\lambda)\big]^{n_1}
f_{A(D_1)}\left(\frac{\lambda+\chi_2(\lambda)}{1-\chi_2(\lambda)}\right)
\left(1-\frac{\chi_2(\lambda)}{1-\chi_2(\lambda)}\,\chi_{1}\left(\frac{\lambda+\chi_2(\lambda)}{1-\chi_2(\lambda)}\right)\right)\\
&=\big[1-\chi_2(\lambda)\big]^{n_1-1}
f_{A(D_1)}\left(\frac{\lambda+\chi_2(\lambda)}{1-\chi_2(\lambda)}\right)
\left(1-\chi_2(\lambda)-\chi_2(\lambda)\,\chi_{1}\left(\frac{\lambda+\chi_2(\lambda)}{1-\chi_2(\lambda)}\right)\right).
\end{align*}
The result now follows from \cref{thm:arccor_char}.
\end{proof}

For the symmetric-arc-corona, the signless Laplacian matrix of the underlying graph of $D_1$ appears in the expression for the characteristic polynomial of $D_1\arcs D_2$ (\cref{thm:arccor_char}(iii)).
When $D_1$ is a symmetric out-regular digraph (i.e., regular graph), we have the following simplification.

\begin{corollary}\label[corollary]{rem:arcb_graph2}
Let $D_1$ and $D_2$ be two digraphs on $n_1$ and $n_2$ vertices, respectively. 
If $D_1$ is a symmetric $r$-out-regular digraph, then 
\[
f_{A(D_1 \arcs D_2)} (\lambda) =\big[f_{A(D_2)}(\lambda)\big]^{rn_1}
\left(1+\chi_{A(D_2)}(\lambda)\right)^{n_1}
f_{A(D_1)}\left(\frac{\lambda-r\,\chi_{A(D_2)}(\lambda)}{1+\chi_{A(D_2)}(\lambda)}\right).
\]
In particular, when $D_1$ is a symmetric out-regular digraph, the adjacency spectrum of $D_1 \arcs D_2$ is completely determined by the adjacency characteristic polynomials of $D_1$ and $D_2$, and the $A$-coronal of $D_2$.
\end{corollary}

\begin{proof}
Since $D_1$ is a symmetric digraph, it has the same adjacency matrix as its underlying graph.
The result follows by \cref{thm:arccor_char} and noting that 
\begin{align*}
\det\left(\lambda I_{n_1}-A(D_1)-\chi_{A(D_2)}(\lambda)\,Q(D_1)\right)
&=
\det\left(\lambda I_{n_1}-A(D_1)-\chi_{A(D_2)}(\lambda)\,(A(D_1)+rI_{n_1})\right)\\
&=
\det\left((\lambda-r\,\chi_{A(D_2)}(\lambda)) I_{n_1}-(1+\chi_{A(D_2)}(\lambda))A(D_1)\right)\\
&=
\left(1+\chi_{A(D_2)}(\lambda)\right)^{n_1}\det\left(\frac{\lambda-r\,\chi_{A(D_2)}(\lambda)}{1+\chi_{A(D_2)}(\lambda)} I_{n_1}-A(D_1)\right)\\
&=\left(1+\chi_{A(D_2)}(\lambda)\right)^{n_1}
f_{A(D_1)}\left(\frac{\lambda-r\,\chi_{A(D_2)}(\lambda)}{1+\chi_{A(D_2)}(\lambda)}\right).
\end{align*}
\end{proof}

\subsection{Laplacian spectrum of arc coronas}
We now give formulas for the Laplacian characteristic polynomials for the three arc-corona constructions.

\begin{theorem}\label[theorem]{thm:arccor_char3}
Let $D_1$ and $D_2$ be two digraphs on $n_1$ and $n_2$ vertices, and $m_1$ and $m_2$ arcs, respectively. 
Suppose $G_1=U(D_1)$ is the underlying graph of $D_1$ with $m_1'$ edges.
Then
\begin{enumerate}[(i)]
    \item $\displaystyle f_{L(D_1 \arcf D_2)} (\lambda)= \big[f_{L(D_2)}(\lambda-1)\big]^{m_1} \det\left(\lambda I_{n_1}-\big(L(D_1)+n_2D_{\rm out}(D_1)\big)-\frac{n_2}{\lambda-1}\,A(D_1)\right)$,        
    \item $\displaystyle f_{L(D_1 \arcb D_2)}(\lambda)= \big[f_{L(D_2)}(\lambda-1)\big]^{m_1} \det\left(\lambda I_{n_1}-\big(L(D_1)+n_2D_{\rm in}(D_1)\big)-\frac{n_2}{\lambda-1}\,A(\overleftarrow{D}_1)\right)$,    
    \item $\displaystyle f_{L(D_1 \arcs D_2)}(\lambda)= \big[f_{L(D_2)}(\lambda-2)\big]^{m_1'}\det\left(\lambda I_{n_1}-\big(L(D_1)+n_2D_{\rm deg}(G_1)\big)-\frac{n_2}{\lambda-2}\,Q(G_1)\right)$.
\end{enumerate}
\end{theorem}

\begin{proof}
The matrices $L(D_1\arcf D_2)$, $L(D_1\arcb D_2)$ and $L(D_1\arcs D_2)$ are permutation similar to
\[
\left[\begin{array}{cc}
L(D_1)+n_2\,D_{\rm out}(D_1)
& -\,\mathbf{1}_{n_2}^T\otimes B_{\rm out}(D_1)\\
-\,\mathbf{1}_{n_2}\otimes B_{\rm in}(D_1)^T
& (L(D_2)+I_{n_2})\otimes I_{m_1}
\end{array}\right],~
\left[\begin{array}{cc}
L(D_1)+n_2\,D_{\rm in}(D_1)
& -\,\mathbf{1}_{n_2}^T\otimes B_{\rm in}(D_1)\\
-\,\mathbf{1}_{n_2}\otimes B_{\rm out}(D_1)^T
& (L(D_2)+I_{n_2})\otimes I_{m_1}
\end{array}\right],
\]
\[
\text{and}~
\left[\begin{array}{cc}
L(D_1)+n_2\,D_{\rm deg}(G_1)
& -\,\mathbf{1}_{n_2}^T\otimes B(G_1)\\
-\,\mathbf{1}_{n_2}\otimes B(G_1)^T
& (L(D_2)+2I_{n_2})\otimes I_{m_1'}
\end{array}\right].
\]
The result now follows from \cref{lm:kronschur} and \cref{cor_addI} and recalling that $\chi_{L(D_2)}(\lambda)=n_2/\lambda$.
\end{proof}

In the case that $D_1$ is an out-regular digraph or a symmetric out-regular digraph, respectively, we have the following two corollaries for the forward-arc corona and symmetric-arc corona.

\begin{corollary}\label[corollary]{lap_arcf1}
Let $D_1$ and $D_2$ be two digraphs on $n_1$ and $n_2$ vertices, respectively. 
If $D_1$ is $r$-out-regular, then 
\[
f_{L(D_1 \arcf D_2)} (\lambda)= \big[f_{L(D_2)}(\lambda-1)\big]^{rn_1} 
\left[\frac{\lambda-n_2-1}{\lambda-1}\right]^{n_1}
f_{L(D_1)}\left(\frac{\lambda^2-(rn_2+1)\lambda}{\lambda-n_2-1}\right).
\]
In particular, when $D_1$ is an out-regular digraph, the Laplacian spectrum of $D_1 \arcf D_2$ is completely determined by the Laplacian characteristic polynomials of $D_1$ and $D_2$.
\end{corollary}

\begin{proof}
Since $m_1=rn_1$ and $A(D_1)=rI_{n_1}-L(D_1)$, by~\Cref{thm:arccor_char3}\,(i), we obtain:
\[
\begin{aligned}
f_{L(D_1 \arcf D_2)}(\lambda)&=\big[f_{L(D_2)}(\lambda-1)\big]^{r n_1}
\det\!\Big(\lambda I_{n_1}-\big(L(D_1)+n_2rI_{n_1}\big)-\tfrac{n_2}{\lambda-1}\,\big(rI_{n_1}-L(D_1)\big)\Big)\\[2mm]
&=\big[f_{L(D_2)}(\lambda-1)\big]^{r n_1}
\det\!\Big({\Big(\lambda-rn_2-\tfrac{rn_2}{\lambda-1}\Big)} I_{n_1}
\;-\;{\Big(1-\tfrac{n_2}{\lambda-1}\Big)}\,L(D_1)\Big).\\
&=\big[f_{L(D_2)}(\lambda-1)\big]^{r n_1}\left[\frac{\lambda-n_2-1}{\lambda-1}\right]^{\!n_1}
\,f_{L(D_1)}\!\left(\frac{\lambda^2-(rn_2+1)\lambda}{\lambda-n_2-1}\right),
\end{aligned}
\]
as required.
\end{proof}

\begin{corollary}\label[corollary]{lap_arcf2}
Let $D_1$ and $D_2$ be two digraphs on $n_1$ and $n_2$ vertices, respectively. 
If $D_1$ is a symmetric $r$-out-regular digraph, then 
\[
f_{L(D_1 \arcs D_2)} (\lambda) =\big[f_{L(D_2)}(\lambda-2)\big]^{rn_1/2}
\left[\frac{\lambda-n_2-2}{\lambda-2}\right]^{n_1}
f_{L(D_1)}\left(\frac{\lambda^2-(rn_2+2)\lambda}{\lambda-n_2-2}\right).
\]
In particular, when $D_1$ is a symmetric out-regular digraph, the Laplacian spectrum of $D_1 \arcs D_2$ is completely determined by the Laplacian characteristic polynomials of $D_1$ and $D_2$.
\end{corollary}

\begin{proof}
Assume $D_1$ is a symmetric $r$-out-regular digraph on $n_1$ vertices.
Then $A(D_1)$ is symmetric and coincides with the adjacency of the underlying
graph $G_1=U(D_1)$, and $m_1'=rn_1/2$, 
\[
L(D_1)=rI_{n_1}-A(D_1),\qquad
D_{\deg}(G_1)=rI_{n_1},\qquad
Q(G_1)=D_{\deg}(G_1)+A(G_1)=rI_{n_1}+A(D_1).
\]
Substituting into \Cref{thm:arccor_char3}(iii), we obtain: 
\[
\begin{aligned}
f_{L(D_1\arcs D_2)}(\lambda)&=\big[f_{L(D_2)}(\lambda-2)\big]^{rn_1/2}\det \big(\lambda I_{n_1}-\big(L(D_1)+n_2 r I_{n_1}\big)-\frac{n_2}{\lambda-2}\,\big(r I_{n_1}+A(D_1)\big)\big)\\[2mm]
&=\big[f_{L(D_2)}(\lambda-2)\big]^{rn_1/2}\det \Bigl({\Bigl(\lambda-n_2 r-\frac{2 n_2 r}{\lambda-2}\Bigr)} I_{n_1}
\;-\;{\Bigl(1-\frac{n_2}{\lambda-2}\Bigr)}\,L(D_1)\Bigr)\\
&=\big[f_{L(D_2)}(\lambda-2)\big]^{rn_1/2}
\left[\frac{\lambda-n_2-2}{\lambda-2}\right]^{n_1}
f_{L(D_1)}\!\left(\frac{\lambda^2-(rn_2+2)\lambda}{\lambda-n_2-2}\right),
\end{aligned}
\]
as claimed.
\end{proof}

\subsection{Signless Laplacian spectrum of arc coronas}
For the signless Laplacian, we have the following formulas for the characteristic polynomials of the three arc-corona constructions.

\begin{theorem}\label[theorem]{thm:arccor_char6}
Let $ D_1 $ and $ D_2 $ be two digraphs on $ n_1 $ and $ n_2 $ vertices, and $ m_1 $ and $ m_2 $ arcs, respectively.
Suppose $G_1=U(D_1)$ is the underlying graph of $D_1$ with $m_1'$ edges.
Then
\begin{enumerate}[(i)]
    \item $\displaystyle f_{Q(D_1 \arcf D_2)} (\lambda)= \big[f_{Q(D_2)}(\lambda-1)\big]^{m_1} \det\left(\lambda I_{n_1}-\big(Q(D_1)+n_2D_{\rm out}(D_1)\big)-\chi_{Q(D_2)}(\lambda-1)\,A(D_1)\right)$,        
    \item $\displaystyle f_{Q(D_1 \arcb D_2)}(\lambda)= \big[f_{Q(D_2)}(\lambda-1)\big]^{m_1} \det\left(\lambda I_{n_1}-\big(Q(D_1)+n_2D_{\rm in}(D_1)\big)-\chi_{Q(D_2)}(\lambda-1)\,A(\overleftarrow{D}_1)\right)$,
    \item $\displaystyle f_{Q(D_1 \arcs D_2)}(\lambda)= \big[f_{Q(D_2)}(\lambda-2)\big]^{m_1'}\det\left(\lambda I_{n_1}-\big(Q(D_1)+n_2D_{\rm deg}(G_1)\big)-\chi_{Q(D_2)}(\lambda-2)\,Q(G_1)\right)$.
\end{enumerate} 
\end{theorem}

\begin{proof}
The matrices $Q(D_1\arcf D_2)$, $Q(D_1\arcb D_2)$ and $Q(D_1\arcs D_2)$ are permutation similar to
\[
\left[\begin{array}{cc}
Q(D_1)+n_2\,D_{\rm out}(D_1)
& \mathbf{1}_{n_2}^T\otimes B_{\rm out}(D_1)\\
\mathbf{1}_{n_2}\otimes B_{\rm in}(D_1)^T
& (Q(D_2)+I_{n_2})\otimes I_{m_1}
\end{array}\right],~
\left[\begin{array}{cc}
Q(D_1)+n_2\,D_{\rm in}(D_1)
& \mathbf{1}_{n_2}^T\otimes B_{\rm in}(D_1)\\
\mathbf{1}_{n_2}\otimes B_{\rm out}(D_1)^T
& (Q(D_2)+I_{n_2})\otimes I_{m_1}
\end{array}\right],
\]
\[
\text{and}~
\left[\begin{array}{cc}
Q(D_1)+n_2\,D_{\rm deg}(G_1)
& \mathbf{1}_{n_2}^T\otimes B(G_1)\\
\mathbf{1}_{n_2}\otimes B(G_1)^T
& (Q(D_2)+2I_{n_2})\otimes I_{m_1'}
\end{array}\right].
\]
The result now follows from \cref{lm:kronschur} and \cref{cor_addI}.
\end{proof}

\begin{corollary}\label[corollary]{slap_arcf1}
Let $D_1$ and $D_2$ be two digraphs on $n_1$ and $n_2$ vertices, and $m_1$ and $m_2$ arcs, respectively. 
If $D_1$ is $r$-out-regular, then   
\[
f_{Q(D_1 \arcf D_2)} (\lambda)=\big[f_{Q(D_2)}(\lambda - 1)\big]^{r n_1} 
\big(1 + \chi_{Q(D_2)}(\lambda - 1)\big)^{n_1} 
f_{Q(D_1)}\!\Bigg(\frac{\lambda-rn_2+r\chi_{Q(D_2)}(\lambda-1)}{1+\chi_{Q(D_2)}(\lambda-1)}\Bigg).
\]
In particular, when $D_1$ is an out-regular digraph, the signless Laplacian spectrum of $D_1 \arcf D_2$ is completely determined by the signless Laplacian characteristic polynomials of $D_1$ and $D_2$ and the $Q$-coronal of $D_2$.
\end{corollary}

\begin{proof}
Assume $D_1$ is $r$-out-regular. Then $D_{\rm out}(D_1)=rI_{n_1}$ and
$Q(D_1)=rI_{n_1}+A(D_1)$.
Moreover $m_1=rn_1$.  
Assume that $\chi\coloneqq \chi_{Q(D_2)}(\lambda-1)$ for brevity. 
Substituting into \Cref{thm:arccor_char6}\,(i), we obtain:
\[
\begin{aligned}
f_{Q(D_1 \arcf D_2)} (\lambda)&=\big[f_{Q(D_2)}(\lambda-1)\big]^{m_1}\!\det\!\Big(\lambda I_{n_1}-\big(Q(D_1)+n_2D_{\rm out}(D_1)\big)-\chi\,A(D_1)\Big)\\
&=\big[f_{Q(D_2)}(\lambda-1)\big]^{rn_1}\!\det\!\Big(\lambda I_{n_1}-\big(rI_{n_1}+A(D_1)+n_2 r I_{n_1}\big)-\chi\,A(D_1)\Big)\\
&=\big[f_{Q(D_2)}(\lambda-1)\big]^{rn_1}\!\det\!\Big(\big(\lambda-r(1+n_2)\big)I_{n_1}-(1+\chi)\,A(D_1)\Big)\\
&=\big[f_{Q(D_2)}(\lambda-1)\big]^{rn_1}\,(1+\chi)^{n_1}\,\!\det\!\Big(\left(\frac{\lambda-r(1+n_2)}{1+\chi}+r\right)I_{n_1}-\left(rI_{n_1}+A(D_1)\right)\Big)\\
&=\big[f_{Q(D_2)}(\lambda-1)\big]^{rn_1}\,(1+\chi)^{n_1}\,
f_{Q(D_1)}\!\left(\frac{\lambda-rn_2+r\chi}{1+\chi}\right),
\end{aligned}
\]
as required.
\end{proof}

\begin{corollary}\label[corollary]{slap_arcf2}
Let $D_1$ and $D_2$ be two digraphs on $n_1$ and $n_2$ vertices, and $m_1$ and $m_2$ arcs, respectively. 
If $D_1$ is a symmetric $r$-out-regular digraph, then  
\[
f_{Q(D_1 \arcs D_2)} (\lambda)=\bigl[f_{Q(D_2)}(\lambda-2)\bigr]^{rn_1/2}
\bigl(1+\chi_{Q(D_2)}(\lambda-2)\bigr)^{n_1}
f_{Q\left(D_1\right)}\left(\frac{\lambda - rn_2}{1+\chi_{Q(D_2)}(\lambda-2)}\right)
\]
In particular, when $D_1$ is a symmetric out-regular digraph, the signless Laplacian spectrum of $D_1 \arcs D_2$ is completely determined by the signless Laplacian characteristic polynomials of $D_1$ and $D_2$ and the $Q$-coronal of $D_2$.
\end{corollary}

\begin{proof}
Let $G_1=U(D_1)$. Since $D_1$ is a symmetric $r$-out-regular digraph, we have  $D_{\deg}(G_1)=rI_{n_1}$, $m_1'=n_1r/2$ and $Q(G_1)=Q(D_1)=rI_{n_1}+A(D_1)$. 
Assume that $\chi\coloneqq \chi_{Q(D_2)}(\lambda-2)$ for brevity. 
Substituting into \Cref{thm:arccor_char6}\,(iii), we obtain:
\begin{align*}
f_{Q(D_1 \arcs D_2)}(\lambda)&
= \big[f_{Q(D_2)}(\lambda-2)\big]^{m_1'}\det\left(\lambda I_{n_1}-\big(Q(D_1)+n_2D_{\rm deg}(G_1)\big)-\chi\,Q(G_1)\right)\\
&= \big[f_{Q(D_2)}(\lambda-2)\big]^{m_1'}\det\left((\lambda-rn_2)I_{n_1}-(1+\chi)Q(D_1)\right)\\
&= \big[f_{Q(D_2)}(\lambda-2)\big]^{m_1'}(1+\chi)^{n_1}\det\left(\frac{\lambda-rn_2}{1+\chi}I_{n_1}-Q(D_1)\right)\\
&= \big[f_{Q(D_2)}(\lambda-2)\big]^{m_1'}(1+\chi)^{n_1}f_{Q\left(D_1\right)}\left(\frac{\lambda - rn_2}{1+\chi}\right),
\end{align*}
as required.
\end{proof}

\section{Further Research}

In this paper, we focus on two specific variants of the corona operation: the vertex and the arc. 
Indeed, there are additional variants that are currently being explored in the literature. 
For instance, \cite{subdivision} studies the subdivision-vertex join,
subdivision-arc join,  
subdivision-vertex neighbourhood corona, and
subdivision-edge neighbourhood corona of two digraphs.



\end{document}